# AN ACCELERATED COMPOSITE GRADIENT METHOD FOR LARGE-SCALE COMPOSITE OBJECTIVE PROBLEMS

MIHAI I. FLOREA* AND SERGIY A. VOROBYOV*

**Abstract.** We introduce a framework, which we denote as the *augmented estimate sequence*, for deriving fast algorithms with provable convergence guarantees. We use this framework to construct a new first-order scheme, the Accelerated Composite Gradient Method (ACGM), for large-scale problems with composite objective structure. ACGM surpasses the state-of-the-art methods for this problem class in terms of provable convergence rate, both in the strongly and non-strongly convex cases, and is endowed with an efficient step size search procedure. We support the effectiveness of our new method with simulation results.

## 1. Introduction.

The field of accelerated first-order methods was created by Nesterov when he introduced his breakthrough Fast Gradient Method (FGM) [10]. FGM was constructed using the simple mathematical machinery of the *estimate sequence* [11] that provides increasingly accurate approximations of a global lower bound tangent at the optimum. Using the estimate sequence, the design process of FGM is straightforward and, by exploiting the problem structure, simultaneously produces state-of-the art convergence guarantees. FGM converges for non-strongly convex objectives at an optimal rate $\mathcal{O}(1/k^2)$ and for strongly convex objectives at a near-optimal rate $\mathcal{O}((1-\sqrt{q})^{-k})$, where $k$ is the iteration index and $q$ is the inverse condition number of the objective [11]. However, FGM requires that the objective be continuously differentiable with Lipschitz gradient, the Lipschitz constant be known in advance, and the problem be unconstrained.

A broad range of problems, including the most common constrained smooth optimization problems, many inverse problems [2], and several classification and reconstruction problems in imaging [6] have a composite structure, i.e., the objective is a sum of a smooth function $f$ with Lipschitz gradient (Lipschitz constant $L_f$) and a simple function $\Psi$, that may embed constraints by including the indicator function of the feasible set. To address the demand for fast algorithms applicable to this problem class, as well as to alleviate the need to know $L_f$ in advance, Nesterov has introduced the Accelerated Multistep Gradient Scheme (AMGS) [13] that relies on *composite gradients* to overcome the limitations of FGM. This algorithm adjusts an estimate of $L_f$ at every step (a process often called "line-search" in the literature [2, 17]) that reflects the local curvature of the function. AMGS is very efficient on sequential machines because the information collected to estimate $L_f$ is reused to advance the algorithm. However, AMGS is slower than FGM when dealing with strongly convex differentiable objectives and the estimation and advancement phases of AMGS's iterations cannot be performed simultaneously on parallel machines.

The Fast Iterative Shrinkage-Thresholding Algorithm (FISTA) [2] solves the parallelization problem by decoupling the advancement phase from the adjustment phase, stalling the former phase only during backtracks. However, FISTA has a fixed $\mathcal{O}(1/k^2)$ provable convergence rate when the objective is strongly convex, and the line-search strategy of the $L_f$ estimate is inferior to that of AMGS in terms of theoretical convergence guarantees. While preparing this manuscript, we became aware of a strongly convex generalization of FISTA, recently introduced in [6], that we designate by FISTA-Chambolle-Pock (FISTA-CP). It has the same convergence guarantees as FGM in both non-strongly and strongly convex cases. The monograph [6] hints at but does not explicitly state any line-search strategy. Furthermore, the utility of both FISTA and FISTA-CP is hampered by the fact that their design process remains obscure and variations are difficult to obtain [5, 20].

*Department of Signal Processing and Acoustics, Aalto University, Espoo, Finland (mihai.florea@aalto.fi, sergiy.vorobyov@aalto.fi).

In this work, we extend the estimate sequence concept to composite objectives by introducing the *augmented estimate sequence*. We formulate a design pattern for first-order methods and we utilize the augmented estimate sequence along with this pattern to derive our Accelerated Composite Gradient Method (ACGM), which has the convergence guarantees of FGM in both the non-strongly and strongly convex cases and is equipped with an explicit adaptive line-search procedure that is decoupled from the advancement phase at every iteration. We further introduce a parallel black-box complexity measure and we show the superiority of our ACGM over the state-of-the-art methods for composite problems in the parallel scenario using theoretical arguments corroborated by simulation results.

**1.1. Assumptions and notation.** We consider the following convex optimization problem

$$\min_{\boldsymbol{x}\in Q} F(\boldsymbol{x}) \overset{\text{def}}{=} f(\boldsymbol{x}) + \Psi(\boldsymbol{x}), \tag{1}$$

where the feasible set $Q \subseteq \mathbb{R}^n$ is closed convex and $\boldsymbol{x}$ is a vector of optimization variables. In this work, we consider only *large-scale* problems [15]. The composite objective $F$ has a non-empty set of optimal points $X^* \subseteq Q$. Function $f : \mathbb{R}^n \to \mathbb{R}$ is convex differentiable with Lipschitz gradient (Lipschitz constant $L_f > 0$) and strong convexity parameter $\mu_f \geq 0$. The regularizer $\Psi : \mathbb{R}^n \to \mathbb{R} \cup \{\infty\}$ is a proper lower semicontinuous convex function with strong convexity parameter $\mu_\Psi$. This implies that $F$ has a strong convexity parameter $\mu = \mu_f + \mu_g$. The regularizer may not be differentiable and it is infinite outside $Q$. However, its proximal map, given by

$$\text{prox}_{\tau\Psi}(\boldsymbol{x}) \overset{\text{def}}{=} \arg\min_{\boldsymbol{z}\in\mathbb{R}^n} \left( \Psi(\boldsymbol{z}) + \frac{1}{2\tau}\|\boldsymbol{z}-\boldsymbol{x}\|_2^2 \right), \quad \boldsymbol{x}\in\mathbb{R}^n, \quad \tau > 0, \tag{2}$$

can be computed with complexity $\mathcal{O}(n)$. The optimization problem is treated by algorithms in a *black-box* setting [9]. The oracle functions are $f(\boldsymbol{x})$, $\nabla f(\boldsymbol{x})$, $\Psi(\boldsymbol{x})$, and $\text{prox}_{\tau\Psi}(\boldsymbol{x})$, with arguments $\boldsymbol{x} \in \mathbb{R}^n$ and $\tau > 0$.

We define a parabola as a quadratic function $\psi : \mathbb{R}^n \to \mathbb{R}$ of the form

$$\psi(\boldsymbol{x}) \overset{\text{def}}{=} p + \frac{\gamma}{2}\|\boldsymbol{x}-\boldsymbol{v}\|_2^2, \quad \boldsymbol{x}\in\mathbb{R}^n, \tag{3}$$

where $\gamma > 0$ gives the curvature, $\boldsymbol{v}$ is the vertex, and $p$ is a constant. We also define $\mathcal{P}$ as the set of all parabolae, $\mathcal{H}$ as the set of all linear functions $h : \mathbb{R}^n \to \mathbb{R}$ (which we denote as hyperplanes), and $\mathcal{G}$ as the set of generalized parabolae, $\mathcal{G} \overset{\text{def}}{=} \mathcal{P} \cup \mathcal{H}$. Parabolae are important in this context because Lipschitz gradient and strong convexity can be defined in terms of parabolic upper bounds and generalized parabolic lower bounds. The set $\mathcal{P}$ is a non-pointed cone which absorbs hyperplanes, that is,

$$\{\alpha_1\phi_1 + \alpha_2\phi_2 + \alpha_3 h \mid \phi_1,\phi_2 \in \mathcal{P},\ h \in \mathcal{H},\ \alpha_1 > 0,\ \alpha_2,\alpha_3 \geq 0\} = \mathcal{P}. \tag{4}$$

The above is a fundamental property in our derivation of ACGM. The proof, provided by Nesterov in [11], is based on double differentiation.

We define two abbreviated expressions, $P_{f,\boldsymbol{y}}(\boldsymbol{x}) \in \mathcal{H}$ and $Q_{f,\gamma,\boldsymbol{y}}(\boldsymbol{x}) \in \mathcal{G}$, as

$$(5) \quad P_{f,\boldsymbol{y}}(\boldsymbol{x}) \overset{\text{def}}{=} f(\boldsymbol{y}) + \langle \nabla f(\boldsymbol{y}), \boldsymbol{x} - \boldsymbol{y} \rangle, \qquad Q_{f,\gamma,\boldsymbol{y}}(\boldsymbol{x}) \overset{\text{def}}{=} P_{f,\boldsymbol{y}}(\boldsymbol{x}) + \frac{\gamma}{2}\|\boldsymbol{x} - \boldsymbol{y}\|_2^2,$$

for any $\boldsymbol{x}, \boldsymbol{y} \in \mathbb{R}^n, \gamma > 0$. Using expression $Q$, we introduce the proximal gradient operator $T_{f,\Psi,L}(\boldsymbol{y})$ as

$$(6) \quad T_{f,\Psi,L}(\boldsymbol{y}) \overset{\text{def}}{=} \arg\min_{\boldsymbol{x} \in \mathbb{R}^n} (Q_{f,L,\boldsymbol{y}}(\boldsymbol{x}) + \Psi(\boldsymbol{x})) = \text{prox}_{\frac{1}{L}\Psi}\left(\boldsymbol{y} - \frac{1}{L}\nabla f(\boldsymbol{y})\right), \ \boldsymbol{y} \in \mathbb{R}^n,$$

where $L > 0$ is a parameter corresponding to the inverse of the step size.

For a given function $\Psi$, we also define the set of composite parabolae

$$(7) \qquad \mathcal{P}_\Psi \overset{\text{def}}{=} \{\psi + c\Psi \mid c \geq 0, \ \psi \in \mathcal{P}\}.$$

**2. Preliminaries.** We present the theoretical building blocks based on concepts from the literature along with a novel complexity measure that we use to construct and analyze the augmented estimate sequence and ACGM.

**2.1. Estimate sequences for composite objectives.** For the class of problems with non-strongly convex objectives, regardless of the optimization algorithm used, the convergence of the iterates can be arbitrarily slow [11] and may only occur in a weak sense [5]. Consequently, we express the convergence rate of first-order schemes on this problem class as the decrease rate of the distance between the objective value and the optimal value. We define a convergence guarantee (provable convergence rate) as the decrease rate of a theoretical upper bound on this distance. When designing algorithms, we index objective values based on iterations. This does not necessarily reflect the actual performance of the algorithm (see subsection 2.3 for discussion). The bound is expressed in terms of points in the domain space (see also [13]) as

$$(8) \qquad A_k(F(\boldsymbol{x}_k) - F(\boldsymbol{x}^*)) \leq \frac{1}{2}\|\boldsymbol{x}_0 - \boldsymbol{x}^*\|_2^2, \quad \boldsymbol{x} \in X^*, \quad k \geq 0,$$

where the weight sequence $\{A_k\}_{k\geq 0}$, $A_k > 0$ for all $k \geq 1$ gives the convergence guarantee. Since the starting point $\boldsymbol{x}_0$ is assumed to be arbitrary, the composite function value $F(\boldsymbol{x}_0)$ may not be finite. Hence, no guarantee can be given for $k = 0$, where $A_0$ is set to 0 to ensure that (8) holds.

The provable convergence rate expression (8) translates to

$$(9) \qquad A_k F(\boldsymbol{x}_k) \leq H_k,$$

where

$$(10) \qquad H_k \overset{\text{def}}{=} A_k F(\boldsymbol{x}^*) + \frac{1}{2}\|\boldsymbol{x}_0 - \boldsymbol{x}^*\|_2^2$$

is the highest allowable upper bound on weighted objective values $A_k F(\boldsymbol{x}_k)$. The convexity of $F$ ensures that there exists a sequence $\{W_k\}_{k\geq 1}$ of global convex lower bounds on $F$, namely

$$(11) \qquad F(\boldsymbol{x}) \geq W_k(\boldsymbol{x}), \quad \boldsymbol{x} \in \mathbb{R}^n, \quad k \geq 1.$$

We define an estimate sequence $\{\psi_k(\boldsymbol{x})\}_{k\geq 0}$ as

$$\psi_k(\boldsymbol{x}) \stackrel{\text{def}}{=} A_k W_k(\boldsymbol{x}) + \frac{\gamma_0}{2}\|\boldsymbol{x}-\boldsymbol{x}_0\|_2^2, \quad 0<\gamma_0\leq 1, \quad k\geq 0, \tag{12}$$

which represents a generalization of the definition given in [13] (where $W_k$ are further required to be linear). Here $\psi_k$, $k \geq 0$, are estimate functions and $\gamma_0$ is the curvature of the initial estimate function $\psi_0$. Since $A_0 = 0$, there is no need to define $W_0$. Both AMGS and FGM are built to maintain the following estimate sequence property[1]

$$A_k F(\boldsymbol{x}_k) \leq \psi_k^*, \quad k \geq 0, \tag{13}$$

where $\psi_k^* \stackrel{\text{def}}{=} \min\limits_{\boldsymbol{x}\in\mathbb{R}^n} \psi_k(\boldsymbol{x})$. The bound (9) follows naturally from (11), (12) and (13). Thus, we have

$$A_k F(\boldsymbol{x}_k) \leq \psi_k^* \leq \psi_k(\boldsymbol{x}^*) \leq H_k, \quad k\geq 0. \tag{14}$$

Estimate sequences were introduced in [13] because $\psi_k^*$ in (13) is easier to maintain than $H_k$ in (9) as an upper bound on the weighted function values $A_k F(\boldsymbol{x}_k)$. However, the gap between $\psi_k^*$ and $H_k$ is large and, as we shall see in subsection 3.2, can be reduced to yield algorithms that are more robust, i.e., applicable without modification to a broader class of problems.

**2.2. Nesterov's accelerated first-order algorithm design pattern.** FGM and AMGS share the structure outlined in Algorithm 1. In the following, we give a detailed explanation of the numerical and functional parameters of this pattern.

**Algorithm 1** A design pattern for Nesterov's first-order accelerated algorithms

$\psi_0(\boldsymbol{x}) = A_0 F(\boldsymbol{x}_0) + \frac{\gamma_0}{2}\|\boldsymbol{x}-\boldsymbol{x}_0\|_2^2$
**for** $k = 0,\ldots,K-1$ **do**
  $L_{k+1} = \mathcal{S}(\boldsymbol{x}_k, \psi_k, A_k, L_k)$
  $a_{k+1} = \mathcal{F}_a(\psi_k, A_k, L_{k+1})$
  $\boldsymbol{y}_{k+1} = \mathcal{F}_y(\boldsymbol{x}_k, \psi_k, A_k, a_{k+1})$
  $A_{k+1} = A_k + a_{k+1}$
  $\boldsymbol{x}_{k+1} = \arg\min\limits_{\boldsymbol{x}\in\mathbb{R}^n} u_{k+1}(\boldsymbol{x})$
  $\psi_{k+1}(\boldsymbol{x}) = \psi_k(\boldsymbol{x}) + a_{k+1} w_{k+1}(\boldsymbol{x})$
**end for**

Algorithm 1 takes as input the starting point $\boldsymbol{x}_0 \in \mathbb{R}^n$, an initial estimate for the Lipschitz constant $L_0 > 0$, the total number of iterations $K > 0$, the initial weight $A_0 \geq 0$, and the curvature $\gamma_0 > 0$. The convergence analysis of AMGS requires that $A_0 = 0$ (also argued in subsection 2.1) and further assumes that $\gamma_0 = 1$. FGM has a different convergence analysis that leads to a slightly modified definition of estimate sequences (for details, see Appendix A). For $A_0 = 0$ and $0 < \gamma_0 \leq 1$, the two estimate sequence definitions coincide. However, FGM assumes that $A_0 > 0$ and, without loss of generality, we set $A_0 = 1$ which replaces the restriction $\gamma_0 \leq 1$ with $\gamma_0 \geq \mu$. The estimate functions of FGM and AMGS take the form

$$\begin{aligned} &\text{FGM:} \quad \psi_k \in \mathcal{P}, \quad \psi_k(\boldsymbol{x}) = p_k + \frac{\gamma_k}{2}\|\boldsymbol{x}-\boldsymbol{v}_k\|_2^2, \quad k \geq 0 \\ &\text{AMGS:} \quad \psi_k \in \mathcal{P}_\Psi, \quad \psi_k(\boldsymbol{x}) = p_k + \frac{1}{2}\|\boldsymbol{x}-\boldsymbol{v}_k\|_2^2 + A_k\Psi(\boldsymbol{x}), \quad k\geq 0. \end{aligned} \tag{15}$$

[1]For FGM, the estimate sequence definition differs slightly, as explained in Appendix A.

At every iteration $k$, the future value of the main iterate $\boldsymbol{x}_{k+1}$ is generated (Algorithm 1, line 7) using majorization minimization, i.e., it is set as the minimum of $u_{k+1}(\boldsymbol{x})$, a local upper bound on $F$. The estimate function is incremented (Algorithm 1, line 8) with a lower bound $w_{k+1}(\boldsymbol{x})$ weighted by $a_{k+1}$. This ensures that the next estimate function $\psi_{k+1}$ adheres to the form found in (15) and in definition (12), where the lower bounds $W_k$ are given by

$$W_k(\boldsymbol{x}) = \frac{1}{A_k}\sum_{i=1}^{k}\left(a_i w_i(\boldsymbol{x})\right), \quad k \geq 1. \tag{16}$$

The line-search procedure $\mathcal{S}$ (Algorithm 1, line 3) outputs an estimate of $L_f$, denoted by $L_{k+1}$ (see section 3 for description). FGM does not use line-search nor the input parameter $L_0$. It assumes that $\Psi(\boldsymbol{x}) = 0$ and that $L_f$ is known in advance. FGM defines the upper bounds based directly on $L_f$.

The test point $\boldsymbol{y}_{k+1}$ and weight $a_{k+1}$ are obtained as functions $\mathcal{F}_a$ and $\mathcal{F}_y$ (Algorithm 1, lines 4 and 5) of the state variables at each iteration. These functions are derived in the algorithm design stage to guarantee that the estimate sequence property (13) carries over to the next iterate, regardless of the algorithmic state. The expressions of functions $\mathcal{F}_a$ and $\mathcal{F}_y$ for FGM and AMGS as well as the lower and upper bounds are listed in Table 1. By replacing the symbols in Table 1 with the corresponding expressions, we recover FGM and AMGS, respectively.

Table 1: Design choices of FGM and AMGS at every iteration $k \geq 0$

| Symbol | In FGM | In AMGS |
|---|---|---|
| $w_{k+1}(\boldsymbol{x})$ | $Q_{f,\mu,\boldsymbol{y}_{k+1}}(\boldsymbol{x})$ | $P_{f,\boldsymbol{x}_{k+1}}(\boldsymbol{x}) + \Psi(\boldsymbol{x})$ |
| $u_{k+1}(\boldsymbol{x})$ | $Q_{f,L_f,\boldsymbol{y}_{k+1}}(\boldsymbol{x})$ | $Q_{f,L_{k+1},\boldsymbol{y}_{k+1}}(\boldsymbol{x}) + \Psi(\boldsymbol{x})$ |
| $\mathcal{F}_a(\psi_k, A_k, L_{k+1})$ | Solution $a > 0$ of $L_f a^2 = (A_k + a)(\gamma_k + \mu a)$ | Solution $a > 0$ of $L_{k+1} a^2 = 2(A_k + a)(1 + \mu A_k)$ |
| $\mathcal{F}_{\boldsymbol{y}}(\boldsymbol{x}_k, \psi_k, A_k, a_{k+1})$ | $\dfrac{A_k\gamma_{k+1}\boldsymbol{x}_k + a_{k+1}\gamma_k\boldsymbol{v}_k}{A_k\gamma_{k+1} + a_{k+1}\gamma_k}$ | $\dfrac{A_k\boldsymbol{x}_k + a_{k+1}\boldsymbol{v}_k}{A_k + a_{k+1}}$ |

**2.3. Wall-clock time units.** When measuring the convergence rate, the prevailing indexing strategies for objective values found in the literature are based on either iterations (e.g., [6, 13, 17, 21]), CPU time in a particular computing environment (e.g., [6, 21]), or the number of calls to a low-level routine that dominates all others in complexity (e.g., [3, 13]). The first approach cannot cope with the diversity of methods studied. For instance, AMGS makes two gradient steps per iteration whereas FISTA makes only one. The second approach does not generalize to the entire problem class while the third does not take into account parallelism.

Recently, with the advent of consumer-grade GPUs and multi-core CPUs, it became possible to run a small number of threads working with large-scale variables on a shared memory machine. With this setup, Nesterov's assumptions on sequential machines still hold: the complexity of computing $f(\boldsymbol{x})$ is comparable to that of $\nabla f(\boldsymbol{x})$ [14]. We denote the amount of wall-clock time required to evaluate $f(\boldsymbol{x})$ or $\nabla f(\boldsymbol{x})$ by 1 wall-clock time unit (WTU). Given that we are dealing with large-scale problems and $\Psi$ is assumed to be simple (e.g. constraints are simple in the most common applications), we attribute a cost of 0 WTU to $\Psi(\boldsymbol{x})$ and $\text{prox}_{\tau\Psi}(\boldsymbol{x})$ calls as well as to individual scalar-vector multiplications and vector additions [15]. Throughout this work, we consider only this *shared memory machine parallel scenario* and thus, we measure algorithm complexity as the objective distance decrease rate when indexed in WTU.

FGM, AMGS, FISTA, and FISTA-CP can be parallelized in black-box form using at most three concurrent computation threads, each thread running on a parallel processing unit, designated as a PPU. A PPU may be itself a single processor or a cluster of processors. However, we assume that all PPUs are identical. Table 2 shows the resource allocation (i.e., at every time unit what computation takes places on which PPU, along with the iteration that computation belongs to) and runtime behavior of the parallelized versions of these algorithms when the line-search, if utilized, is successful. When the search parameters are tuned properly, backtracks rarely occur, making this situation the most common one for the algorithms presented. Parallelization involves the technique of speculative execution [8] whereby the validation phase of the search takes place in parallel with the advancement phase of the next iteration. Only FISTA is parallelizable on the function call level. When no backtracks occur, it has a running time complexity of 1 WTU per iteration, along with FGM and FISTA-CP, which lack an explicit search scheme. On the other hand, AMGS requires 2 WTU per iteration in this case.

Table 2: Resource allocation and runtime behavior of parallel black-box FGM, FISTA-CP, AMGS and FISTA when no backtracks occur (iteration $k \geq 1$ starts at time T)

| Method | WTU | PPU 1 | | PPU 2 | | PPU 3 | |
|---|---|---|---|---|---|---|---|
| | | Comp. | Iter. | Comp. | Iter. | Comp. | Iter. |
| FGM | T | $\nabla f(\boldsymbol{y}_{k+1})$ | k | Idle | | Idle | |
| | T + 1 | $\nabla f(\boldsymbol{y}_{k+2})$ | k + 1 | Idle | | Idle | |
| FISTA-CP | T | $\nabla f(\boldsymbol{y}_{k+1})$ | k | Idle | | Idle | |
| | T + 1 | $\nabla f(\boldsymbol{y}_{k+2})$ | k + 1 | Idle | | Idle | |
| AMGS | T | $\nabla f(\boldsymbol{y}_{k+1})$ | k | Idle | | Idle | |
| | T + 1 | $\nabla f(\boldsymbol{x}_{k+1})$ | k | Idle | | Idle | |
| | T + 2 | $\nabla f(\boldsymbol{y}_{k+2})$ | k + 1 | Idle | | Idle | |
| FISTA | T | $\nabla f(\boldsymbol{y}_{k+1})$ | k | $f(\boldsymbol{y}_{k+1})$ | k | $f(\boldsymbol{x}_k)$ | k - 1 |
| | T + 1 | $\nabla f(\boldsymbol{y}_{k+2})$ | k + 1 | $f(\boldsymbol{y}_{k+2})$ | k + 1 | $f(\boldsymbol{x}_{k+1})$ | k |
| | T + 2 | $\nabla f(\boldsymbol{y}_{k+3})$ | k + 2 | $f(\boldsymbol{y}_{k+3})$ | k + 2 | $f(\boldsymbol{x}_{k+2})$ | k + 1 |

When a backtrack occurs, function and gradient values of points that change have to be recomputed, stalling the entire multi-threaded system accordingly. It follows that additional backtracks have the same cost. A backtrack costs 2 WTU for AMGS and only 1 WTU for FISTA (see Table 3).

**2.4. Composite gradient.** FGM has a provable convergence rate, measured in the previously introduced WTU, that is superior to that of AMGS (for details and proof, see Appendix B), but it requires that the objective be differentiable. A link between FGM and AMGS has been provided in [13] by means of the *composite gradient*, defined as

$$g_L(\boldsymbol{y}) \stackrel{\text{def}}{=} L\,(\boldsymbol{y} - T_{f,\Psi,L}(\boldsymbol{y})), \quad \boldsymbol{y} \in \mathbb{R}^n, \quad L > 0. \tag{17}$$

As we shall see in (30), there is no need specify functional parameters. The composite gradient substitutes the gradient for composite functions and shares many of its

Table 3: Resource allocation and runtime behavior of parallel black-box AMGS and FISTA when a single backtrack occurs (iteration $k \geq 1$ starts at time T)

| Method | WTU | PPU 1 | | PPU 2 | | PPU 3 | |
|---|---|---|---|---|---|---|---|
| | | Comp. | Iter. | Comp. | Iter. | Comp. | Iter. |
| AMGS | T | $\nabla f(\boldsymbol{y}_{k+1})$ | k | Idle | | Idle | |
| | T + 1 | $\nabla f(\boldsymbol{x}_{k+1})$ | k | Idle | | Idle | |
| | T + 2 | $\nabla f(\boldsymbol{y}_{k+1})$ | k | Idle | | Idle | |
| | T + 3 | $\nabla f(\boldsymbol{x}_{k+1})$ | k | Idle | | Idle | |
| | T + 4 | $\nabla f(\boldsymbol{y}_{k+2})$ | k + 1 | Idle | | Idle | |
| FISTA | T | $\nabla f(\boldsymbol{y}_{k+1})$ | k | $f(\boldsymbol{y}_{k+1})$ | k | $f(\boldsymbol{x}_k)$ | k - 1 |
| | T + 1 | $\nabla f(\boldsymbol{y}_{k+2})$ | k + 1 | $f(\boldsymbol{y}_{k+2})$ | k + 1 | $f(\boldsymbol{x}_{k+1})$ | k |
| | T + 2 | $\nabla f(\boldsymbol{y}_{k+2})$ | k + 1 | $f(\boldsymbol{y}_{k+2})$ | k + 1 | $f(\boldsymbol{x}_{k+1})$ | k |
| | T + 3 | $\nabla f(\boldsymbol{y}_{k+3})$ | k + 2 | $f(\boldsymbol{y}_{k+3})$ | k + 2 | $f(\boldsymbol{x}_{k+2})$ | k + 1 |

properties. Most notably, the descent update (Algorithm 1, line 7) in FGM, given by

$$\boldsymbol{x}_{k+1} = \boldsymbol{y}_{k+1} - \frac{1}{L_f}\nabla f(\boldsymbol{y}_{k+1}), \tag{18}$$

can be written similarly in AMGS using the composite gradient as

$$\boldsymbol{x}_{k+1} = \boldsymbol{y}_{k+1} - \frac{1}{L_{k+1}} g_{L_{k+1}}(\boldsymbol{y}_{k+1}). \tag{19}$$

In addition, the descent rule [11], which for FGM takes the form of

$$f(\boldsymbol{x}_{k+1}) \leq f(\boldsymbol{y}_{k+1}) - \frac{1}{2L_f}\|\nabla f(\boldsymbol{y}_{k+1})\|_2^2, \tag{20}$$

is obeyed by the composite gradient in AMGS as well (see Lemma 2), that is,

$$F(\boldsymbol{x}_{k+1}) \leq F(\boldsymbol{y}_{k+1}) - \frac{1}{2L_{k+1}}\|g_{L_{k+1}}(\boldsymbol{y}_{k+1})\|_2^2. \tag{21}$$

These properties suggest that FGM could be applied to composite objectives simply by replacing the gradient call with a composite gradient call, yielding an algorithm that has the superior convergence guarantees of FGM and the applicability of AMGS.

**3. ACGM.** The convergence analysis of FGM in [11] requires only two properties of the gradient to hold: the descent rule (20) and the supporting generalized parabola condition, i.e., $Q_{f,\mu,\boldsymbol{y}_{k+1}}$ is a lower bound on function $f$ for all $k \geq 0$. However, the extension of $Q_{f,\mu,\boldsymbol{y}_{k+1}}(\boldsymbol{x})$ to composite gradients, written as

$$F(\boldsymbol{y}_{k+1}) + \langle g_{L_{k+1}}(\boldsymbol{y}_{k+1}), \boldsymbol{x} - \boldsymbol{y}_{k+1}\rangle + \frac{\mu}{2}\|\boldsymbol{x} - \boldsymbol{y}_{k+1}\|_2^2, \tag{22}$$

is *not guaranteed* to be a valid lower bound on $F$ for *any* value of $L_{k+1} > 0$. Hence, this convergence analysis of FGM does not apply to composite objectives. In the following subsection, we seek a suitable replacement for the FGM supporting generalized parabolae, bearing in mind that the accuracy of the lower bounds at every iteration positively impacts the convergence rate of the algorithm (also argued in Appendix B).

**3.1. Relaxed lower bound.** At every iteration $k$, the lower bound in FGM takes the form of an approximate second order Taylor expansion of $f$ at $\boldsymbol{y}_{k+1}$. For ACGM, we produce a similar lower bound on $F$ by transferring all strong convexity from $\Psi$ to $f$ as

$$f'(\boldsymbol{x}) \stackrel{\text{def}}{=} f(\boldsymbol{x}) + \frac{\mu_\Psi}{2}\|\boldsymbol{x} - \boldsymbol{x}_0\|_2^2, \qquad \Psi'(\boldsymbol{x}) \stackrel{\text{def}}{=} \Psi(\boldsymbol{x}) - \frac{\mu_\Psi}{2}\|\boldsymbol{x} - \boldsymbol{x}_0\|_2^2. \tag{23}$$

The center of strong convexity in (23) can be any point in $\mathbb{R}^n$. We choose $\boldsymbol{x}_0$ only for convenience. Function $f'$ has Lipschitz gradient with constant $L_{f'} = L_f + \mu_\Psi$ and strong convexity parameter $\mu_{f'} = \mu$. Naturally, this transfer does not alter the objective function

$$F(\boldsymbol{x}) = f(\boldsymbol{x}) + \Psi(\boldsymbol{x}) = f'(\boldsymbol{x}) + \Psi'(\boldsymbol{x}) \tag{24}$$

and gives rise to the following remarkable property.

Proposition 1. *By transferring convexity as in (23) we have*

$$Q_{f',L+\mu_\Psi,\boldsymbol{y}}(\boldsymbol{x}) = Q_{f,L,\boldsymbol{y}}(\boldsymbol{x}) + \frac{\mu_\Psi}{2}\|\boldsymbol{x} - \boldsymbol{x_0}\|_2^2, \quad \boldsymbol{x}, \boldsymbol{y} \in \mathbb{R}^n, \quad L > 0. \tag{25}$$

*Proof.* By expanding $Q_{f',L+\mu_\Psi,\boldsymbol{y}}$ we obtain

$$\begin{aligned} Q_{f',L+\mu_\Psi,\boldsymbol{y}}(\boldsymbol{x}) &= f'(\boldsymbol{y}) + \langle \nabla f'(\boldsymbol{y}), \boldsymbol{x} - \boldsymbol{y}\rangle + \frac{L+\mu_\Psi}{2}\|\boldsymbol{x} - \boldsymbol{y}\|_2^2 \\ &= f(\boldsymbol{y}) + \frac{\mu_\Psi}{2}\|\boldsymbol{y} - \boldsymbol{x}_0\|_2^2 + \langle \nabla f(\boldsymbol{y}) + \mu_\Psi(\boldsymbol{y} - \boldsymbol{x}_0), \boldsymbol{x} - \boldsymbol{y}\rangle + \\ &\quad + \frac{L+\mu_\Psi}{2}\|\boldsymbol{x} - \boldsymbol{y}\|_2^2 - \frac{\mu_\Psi}{2}\|\boldsymbol{x} - \boldsymbol{x}_0\|_2^2 + \frac{\mu_\Psi}{2}\|\boldsymbol{x} - \boldsymbol{x}_0\|_2^2 \\ &= f(\boldsymbol{y}) + \langle \nabla f(\boldsymbol{y}), \boldsymbol{x} - \boldsymbol{y}\rangle + \frac{L}{2}\|\boldsymbol{x} - \boldsymbol{y}\|_2^2 + \frac{\mu_\Psi}{2}\|\boldsymbol{x} - \boldsymbol{x}_0\|_2^2, \end{aligned} \tag{26}$$

Rewriting (26) using the definition of $Q$ in (5) completes the proof. ∎

From Proposition 1 it follows that, by setting $L'_{k+1} \stackrel{\text{def}}{=} L_{k+1} + \mu_\Psi$, the descent conditions for $f$ and $f'$ at every iteration $k$ are equivalent, namely

$$f(\boldsymbol{x}_{k+1}) \leq Q_{f,L_{k+1},\boldsymbol{y}_{k+1}}(\boldsymbol{x}_{k+1}) \quad \Leftrightarrow \quad f'(\boldsymbol{x}_{k+1}) \leq Q_{f',L'_{k+1},\boldsymbol{y}_{k+1}}(\boldsymbol{x}_{k+1}). \tag{27}$$

When designing ACGM, since we assume no upper bound on $\Psi$, we have to choose a composite parabolic upper bound on $F$ at every iteration $k \geq 0$, that is,

$$u_{k+1}(\boldsymbol{x}) = Q_{f,L_{k+1},\boldsymbol{y}_{k+1}}(\boldsymbol{x}) + \Psi(\boldsymbol{x}), \quad \boldsymbol{x} \in \mathbb{R}^n. \tag{28}$$

From Proposition 1 we can also see that the strong convexity transfer in (25) does not alter the upper bound, namely

$$u_{k+1}(\boldsymbol{x}) = Q_{f,L_{k+1},\boldsymbol{y}_{k+1}}(\boldsymbol{x}) + \Psi(\boldsymbol{x}) = Q_{f',L'_{k+1},\boldsymbol{y}_{k+1}}(\boldsymbol{x}) + \Psi'(\boldsymbol{x}), \quad \boldsymbol{x} \in \mathbb{R}^n. \tag{29}$$

Hence, the update in line 7 in Algorithm 1 remains unchanged as well,

$$\boldsymbol{x}_{k+1} = T_{f,\Psi,L_{k+1}}(\boldsymbol{y}_{k+1}) = T_{f',\Psi',L'_{k+1}}(\boldsymbol{y}_{k+1}). \tag{30}$$

LEMMA 2. *If the descent condition (27) holds at iteration $k \geq 0$, then the objective $F$ is lower bounded as*

$$F(\boldsymbol{x}) \geq \mathcal{R}_{L'_{k+1},\boldsymbol{y}_{k+1}}(\boldsymbol{x}), \quad \boldsymbol{x} \in \mathbb{R}^n, \tag{31}$$

*where we denote with $\mathcal{R}_{L'_{k+1},\boldsymbol{y}_{k+1}}(\boldsymbol{x})$ the* relaxed supporting generalized parabola *of $F$ at $\boldsymbol{y}_{k+1}$ using inverse step size $L'_{k+1}$,*

$$\begin{aligned} \mathcal{R}_{L'_{k+1},\boldsymbol{y}_{k+1}}(\boldsymbol{x}) \overset{\text{def}}{=} &\ F(\boldsymbol{x}_{k+1}) + \frac{1}{2L'_{k+1}}\|g_{L'_{k+1}}(\boldsymbol{y}_{k+1})\|_2^2 + \\ &+ \langle g_{L'_{k+1}}(\boldsymbol{y}_{k+1}), \boldsymbol{x} - \boldsymbol{y}_{k+1}\rangle + \frac{\mu}{2}\|\boldsymbol{x} - \boldsymbol{y}_{k+1}\|_2^2, \quad \boldsymbol{x} \in \mathbb{R}^n, \end{aligned} \tag{32}$$

*with $\boldsymbol{x}_{k+1}$ given by (30).*

*Proof.* From the strong convexity of $f'$, we have a supporting generalized parabola at $\boldsymbol{y}_{k+1}$, given by

$$f'(\boldsymbol{x}) \geq f'(\boldsymbol{y}_{k+1}) + \langle \nabla f'(\boldsymbol{y}_{k+1}), \boldsymbol{x} - \boldsymbol{y}_{k+1}\rangle + \frac{\mu}{2}\|\boldsymbol{x} - \boldsymbol{y}_{k+1}\|_2^2, \quad \boldsymbol{x} \in \mathbb{R}^n. \tag{33}$$

The first-order optimality condition of (6) implies that there exists a subgradient $\boldsymbol{\xi}$ of function $\Psi'$ at point $\boldsymbol{x}_{k+1}$ such that

$$g_{L'_{k+1}}(\boldsymbol{y}_{k+1}) = \nabla f'(\boldsymbol{y}_{k+1}) + \boldsymbol{\xi}. \tag{34}$$

From the convexity of $\Psi'$, we have a supporting hyperplane at $\boldsymbol{x}_{k+1}$, which satisfies

$$\begin{aligned} \Psi'(\boldsymbol{x}) \geq &\ \Psi'(\boldsymbol{x}_{k+1}) + \langle \boldsymbol{\xi}, \boldsymbol{x} - \boldsymbol{x}_{k+1}\rangle \\ = &\ \Psi'(\boldsymbol{x}_{k+1}) + \langle g_{L'_{k+1}}(\boldsymbol{y}_{k+1}) - \nabla f'(\boldsymbol{y}_{k+1}), \boldsymbol{x} - \boldsymbol{x}_{k+1}\rangle, \quad \boldsymbol{x} \in \mathbb{R}^n. \end{aligned} \tag{35}$$

By adding together (33), (35), and the descent condition for $f'$ (27), we obtain the desired result (31). ∎

The relaxed supporting generalized parabola uses recent information ($\boldsymbol{y}_{k+1}$ and $\boldsymbol{x}_{k+1}$) and has the same curvature as the FGM lower bounds. However, the estimate sequence property cannot be maintained across iterations using the pattern in Algorithm 1 with such a loose global lower bound on $F$.

**3.2. Augmented estimate sequence.** Recall that the estimate sequence property (13) produces a large gap between $\psi^*$ and $H_k$, resulting in algorithms that are excessively stringent. To remedy this problem, we introduce the *augmented estimate sequence* $\{\psi'_k(\boldsymbol{x})\}_{k\geq 0}$ with a gap low enough to accommodate relaxed supporting generalized parabolae as lower bounds. Using the notation and conventions from subsection 2.1, we define the augmented estimate sequence as

$$\psi'_k(\boldsymbol{x}) \overset{\text{def}}{=} \psi_k(\boldsymbol{x}) + A_k(F(\boldsymbol{x}^*) - W_k(\boldsymbol{x}^*)), \quad k \geq 0. \tag{36}$$

Augmentation consists only of adding a non-negative constant (due to the lower bound property of $W_k$) to the estimate function, thus preserving its curvature and vertex. The augmented estimate sequence property, given as

$$A_k F(\boldsymbol{x}_k) \leq \psi'^*_k, \quad k \geq 0, \tag{37}$$

can be used to derive the provable convergence rate because, along with definitions (10), (12), and (36), it implies that

$$A_k F(\boldsymbol{x}_k) \leq \psi'^*_k = \psi^*_k + A_k(F(\boldsymbol{x}^*) - W_k(\boldsymbol{x}^*)) \leq \psi^*_k + H_k - \psi_k(\boldsymbol{x}^*) \leq H_k. \tag{38}$$

**3.3. Formulating ACGM.** We proceed with the design of our method, ACGM, based on the pattern presented in Algorithm 1. The building blocks are:

1. The augmented estimate sequence property (37) .
2. The composite parabolic upper bounds (28). This implies that line 7 in Algorithm 1 is the proximal gradient step (30).
3. The relaxed supporting generalized parabola lower bounds from Lemma 2,

$$w_{k+1}(\boldsymbol{x}) = \mathcal{R}_{L'_{k+1},\boldsymbol{y}_{k+1}}(\boldsymbol{x}), \quad \boldsymbol{x} \in \mathbb{R}^n, \quad k \geq 0. \tag{39}$$

For the relaxed supporting generalized parabola to be a valid global lower bound on $F$, Lemma 2 requires that, at every iteration $k$, the descent condition for $f$ (27) holds. This is assured in the worst case when $L_{k+1} \geq L_f$.

We construct ACGM by induction. First, we presume that at an arbitrary iteration $k \geq 0$, the augmented estimate sequence property is satisfied. Then, we devise update rules for $a_{k+1}$ and $\boldsymbol{y}_{k+1}$ such that the augmented estimate sequence property is guaranteed to hold at iteration $k+1$ as well, for any algorithmic state. Given that initially, $A_0 F(\boldsymbol{x}_0) - \psi_0'^* = 0$, a sufficient condition for this guarantee is that the gap between the weighted function values and the augmented estimate sequence optimum is monotonically decreasing, namely

$$A_{k+1}F(\boldsymbol{x}_{k+1}) - \psi_{k+1}'^* \leq A_k F(\boldsymbol{x}_k) - \psi_k'^*, \quad k \geq 0. \tag{40}$$

Since the initial estimate function is a parabola and the lower bounds are generalized parabolae, we can write the estimate function at any iteration $k$, along with its augmented variant, as the following parabolae:

$$\psi_k(\boldsymbol{x}) = \psi_k^* + \frac{\gamma_k}{2}\|\boldsymbol{x} - \boldsymbol{v_k}\|_2^2, \qquad \psi_k'(\boldsymbol{x}) = \psi_k'^* + \frac{\gamma_k}{2}\|\boldsymbol{x} - \boldsymbol{v_k}\|_2^2. \tag{41}$$

The gap between $A_k F(\boldsymbol{x}_k)$ and $\psi_k'^*$ can be expressed as

$$\begin{aligned} A_k F(\boldsymbol{x}_k) - \psi_k'^* &\overset{(36)}{=} A_k(F(\boldsymbol{x}_k) - F(\boldsymbol{x}^*)) + A_k W_k(\boldsymbol{x}^*) - \psi_k^* \\ &\overset{(12)}{=} A_k(F(\boldsymbol{x}_k) - F(\boldsymbol{x}^*)) + \psi_k(\boldsymbol{x}^*) - \psi_k^* - \frac{\gamma_0}{2}\|\boldsymbol{x}^* - \boldsymbol{x}_0\|_2^2 \\ &\overset{(41)}{=} A_k(F(\boldsymbol{x}_k) - F(\boldsymbol{x}^*)) + \frac{\gamma_k}{2}\|\boldsymbol{v}_k - \boldsymbol{x}^*\|_2^2 - \frac{\gamma_0}{2}\|\boldsymbol{x}^* - \boldsymbol{x}_0\|_2^2, \quad k \geq 0. \end{aligned} \tag{42}$$

We define the gap sequence $\{\Delta_k\}_{k\geq 0}$ as

$$\Delta_k \overset{\text{def}}{=} A_k(F(\boldsymbol{x}_k) - F(\boldsymbol{x}^*)) + \frac{\gamma_k}{2}\|\boldsymbol{v}_k - \boldsymbol{x}^*\|_2^2. \tag{43}$$

With the quantity $\frac{\gamma_0}{2}\|\boldsymbol{x}^* - \boldsymbol{x}_0\|_2^2$ being constant across iterations, the sufficient condition (40) can be written as

$$\Delta_{k+1} \leq \Delta_k, \quad k \geq 0. \tag{44}$$

THEOREM 3. *If at iteration $k \geq 0$, the descent condition for $f$ (27) holds, then*

$$\Delta_{k+1} + \mathcal{A}_k + \mathcal{B}_k \leq \Delta_k, \tag{45}$$

*where subexpressions $\mathcal{A}_k$, $\mathcal{B}_k$, and the reduced composite gradient $\boldsymbol{G}_k$ are defined as*

$$\mathcal{A}_k \overset{\text{def}}{=} \left(\frac{A_{k+1}}{2L'_{k+1}} - \frac{a_{k+1}^2}{2\gamma_{k+1}}\right)\|g_{L'_{k+1}}(\boldsymbol{y}_{k+1})\|_2^2, \tag{46}$$

$$\mathcal{B}_k \overset{\text{def}}{=} \frac{1}{\gamma_{k+1}}\langle \boldsymbol{G}_k, A_k\gamma_{k+1}\boldsymbol{x}_k + a_{k+1}\gamma_k\boldsymbol{v}_k - (A_k\gamma_{k+1} + a_{k+1}\gamma_k)\boldsymbol{y}_{k+1}\rangle, \tag{47}$$

$$\boldsymbol{G}_k \overset{\text{def}}{=} g_{L'_{k+1}}(\boldsymbol{y}_{k+1}) - \mu\boldsymbol{y}_{k+1}. \tag{48}$$

*Proof.* By applying Lemma 2 at iteration $k$ using $\boldsymbol{x}_k$ and $\boldsymbol{x}^*$ as values of $\boldsymbol{x}$, we obtain

$$F(\boldsymbol{x}_k) - F(\boldsymbol{x}_{k+1}) \geq \frac{1}{2L'_{k+1}}\|g_{L'_{k+1}}(\boldsymbol{y}_{k+1})\|_2^2 + \\ + \langle g_{L'_{k+1}}(\boldsymbol{y}_{k+1}), \boldsymbol{x}_k - \boldsymbol{y}_{k+1}\rangle + \frac{\mu}{2}\|\boldsymbol{x}_k - \boldsymbol{y}_{k+1}\|_2^2, \tag{49}$$

$$F(\boldsymbol{x}^*) - F(\boldsymbol{x}_{k+1}) \geq \frac{1}{2L'_{k+1}}\|g_{L'_{k+1}}(\boldsymbol{y}_{k+1})\|_2^2 + \\ + \langle g_{L'_{k+1}}(\boldsymbol{y}_{k+1}), \boldsymbol{x}^* - \boldsymbol{y}_{k+1}\rangle + \frac{\mu}{2}\|\boldsymbol{x}^* - \boldsymbol{y}_{k+1}\|_2^2. \tag{50}$$

Combining these two inequalities as $A_k \cdot$ (49) $+ a_{k+1} \cdot$ (50), we get

$$A_k(F(\boldsymbol{x}_k) - F(\boldsymbol{x}^*)) - A_{k+1}(F(\boldsymbol{x}_{k+1}) - F(\boldsymbol{x}^*)) \geq \mathcal{C}_k, \tag{51}$$

where the lower bound $\mathcal{C}_k$ is defined as

$$\mathcal{C}_k \overset{\text{def}}{=} \frac{A_{k+1}}{2L'_{k+1}}\|g_{L'_{k+1}}(\boldsymbol{y}_{k+1})\|_2^2 + \frac{\mu}{2}A_k\|\boldsymbol{x}_k - \boldsymbol{y}_{k+1}\|_2^2 + \frac{\mu}{2}a_{k+1}\|\boldsymbol{x}^* - \boldsymbol{y}_{k+1}\|_2^2 + \\ + \langle g_{L'_{k+1}}(\boldsymbol{y}_{k+1}), A_k(\boldsymbol{x}_k - \boldsymbol{y}_{k+1}) + a_{k+1}(\boldsymbol{x}^* - \boldsymbol{y}_{k+1})\rangle. \tag{52}$$

By substituting the relaxed supporting generalized parabola lower bound (39) in the estimate sequence update (Algorithm 1, line 8), we obtain (e.g. by successive derivation) the recursion rules for curvature and vertices as

$$\gamma_{k+1} = \gamma_k + a_{k+1}\mu, \tag{53}$$

$$\gamma_{k+1}\boldsymbol{v}_{k+1} = \gamma_k\boldsymbol{v}_k - a_{k+1}(g_{L'_{k+1}}(\boldsymbol{y}_{k+1}) - \mu\boldsymbol{y}_{k+1}) = \gamma_k\boldsymbol{v}_k - a_{k+1}\boldsymbol{G}_k. \tag{54}$$

Using the reduced composite gradient, expression $\mathcal{C}_k$ becomes

$$\begin{aligned}\mathcal{C}_k &= \mathcal{A}_k + \mathcal{D}_k + \frac{\|a_{k+1}\boldsymbol{G}_k\|_2^2}{2\gamma_{k+1}} + \left(\frac{a_{k+1}\mu}{\gamma_{k+1}} - A_{k+1}\right)\langle \boldsymbol{G}_k, \boldsymbol{y}_{k+1}\rangle + A_k\langle \boldsymbol{G}_k, \boldsymbol{x}_k\rangle \\ &= \mathcal{A}_k + \mathcal{D}_k + \frac{1}{\gamma_{k+1}}\langle \boldsymbol{G}_k, \frac{a_{k+1}^2}{2}\boldsymbol{G}_k + A_k\gamma_{k+1}\boldsymbol{x}_k - (a_{k+1}\gamma_k + A_k\gamma_{k+1})\boldsymbol{y}_{k+1}\rangle \\ &= \mathcal{A}_k + \mathcal{B}_k + \mathcal{D}_k + \mathcal{E}_k, \end{aligned} \tag{55}$$

where, for brevity, we define subexpressions $\mathcal{D}_k$ and $\mathcal{E}_k$ as

$$\begin{aligned}\mathcal{D}_k &\overset{\text{def}}{=} \frac{\mu a_{k+1}}{2}\|\boldsymbol{x}^*\|_2^2 + a_{k+1}\langle \boldsymbol{G}_k, \boldsymbol{x}^*\rangle + \frac{\mu A_k}{2}\|\boldsymbol{x}_k\|_2^2 + \frac{a_{k+1}^2\mu^2}{2\gamma_{k+1}}\|\boldsymbol{y}_{k+1}\|_2^2, \\ \mathcal{E}_k &\overset{\text{def}}{=} \frac{1}{2\gamma_{k+1}}\langle a_{k+1}\boldsymbol{G}_k, a_{k+1}\boldsymbol{G}_k - 2\gamma_k\boldsymbol{v}_k\rangle.\end{aligned} \tag{56}$$

Subexpression $\mathcal{E}_k$ can be lower bounded as

$$\begin{aligned}\mathcal{E}_k &\overset{(54)}{=} \frac{1}{2\gamma_{k+1}}\langle \gamma_k \boldsymbol{v}_k - \gamma_{k+1}\boldsymbol{v}_{k+1}, \boldsymbol{v}_k - \gamma_{k+1}\boldsymbol{v}_{k+1} - 2\gamma_k \boldsymbol{v}_k\rangle \\ &= \frac{\gamma_{k+1}}{2}\|\boldsymbol{v}_{k+1}\|_2^2 - \frac{\gamma_k^2}{2\gamma_{k+1}}\|\boldsymbol{v}_k\|_2^2 \\ &\overset{(53)}{\geq} \frac{\gamma_{k+1}}{2}\|\boldsymbol{v}_{k+1}\|_2^2 - \frac{\gamma_k}{2}\|\boldsymbol{v}_k\|_2^2. \end{aligned} \tag{57}$$

Finally, using (57) in (51), we obtain that

$$\begin{aligned} &A_k(F(\boldsymbol{x}_k) - F(\boldsymbol{x}^*)) - A_{k+1}(F(\boldsymbol{x}_{k+1}) - F(\boldsymbol{x}^*)) \geq \\ &\geq \frac{\gamma_{k+1}}{2}\|\boldsymbol{v}_{k+1} - \boldsymbol{x}^*\|_2^2 - \frac{\gamma_k}{2}\|\boldsymbol{v}_k - \boldsymbol{x}^*\|_2^2 + \mathcal{A}_k + \mathcal{B}_k + \\ &\quad + \frac{\mu A_k}{2}\|\boldsymbol{x}_k\|_2^2 + \frac{a_{k+1}^2\mu^2}{2\gamma_{k+1}}\|\boldsymbol{y}_{k+1}\|_2^2. \end{aligned} \tag{58}$$

The square (non-negative) terms $\|\boldsymbol{x}_k\|_2^2$ and $\|\boldsymbol{y}_{k+1}\|_2^2$ can both be arbitrarily small and can be left out. Rearranging terms completes the proof. □

Theorem 3 implies that the sufficient condition (44) holds if, for any algorithmic state, $\mathcal{A}_k \geq 0$ and $\mathcal{B}_k \geq 0$. The simplest way to achieve this is by maintaining the following conditions:

$$(A_k\gamma_{k+1} + a_{k+1}\gamma_k)\boldsymbol{y}_{k+1} = A_k\gamma_{k+1}\boldsymbol{x}_k + a_{k+1}\gamma_k\boldsymbol{v}_k, \tag{59}$$

$$A_{k+1}\gamma_{k+1} \geq L'_{k+1}a_{k+1}^2 = (L_{k+1} + \mu_\Psi)a_{k+1}^2. \tag{60}$$

Therefore, function $\mathcal{F}_y$ is given by

$$\boldsymbol{y}_{k+1} = \mathcal{F}_y(\boldsymbol{x_k}, \psi_k, A_k, a_{k+1}) = \frac{A_k\gamma_{k+1}\boldsymbol{x}_k + a_{k+1}\gamma_k\boldsymbol{v}_k}{A_k\gamma_{k+1} + a_{k+1}\gamma_k}, \tag{61}$$

where $\gamma_{k+1}$ is obtained from (53).

We choose the most aggressive accumulated weight update by enforcing equality in (60) and ensuring that $\gamma_{k+1}$ is as large as possible by setting $\gamma_0 = 1$. Update (60) becomes

$$(L_{k+1} + \mu_\Psi)a_{k+1}^2 = A_{k+1}\gamma_{k+1} \overset{(53)}{=} (A_k + a_{k+1})(\gamma_k + \mu a_{k+1}). \tag{62}$$

Given that $a_{k+1}, L_{k+1} > 0$ and $A_k \geq 0$, we can write $\mathcal{F}_a$ in closed form as

$$a_{k+1} = \mathcal{F}_a(\psi_k, A_k, L_{k+1}) = \frac{\gamma_k + A_k\mu + \sqrt{(\gamma_k + A_k\mu)^2 + 4(L_{k+1} - \mu_f)A_k\gamma_k}}{2(L_{k+1} - \mu_f)} \tag{63}$$

There is no need to compute the composite gradient or the reduced composite gradient explicitly. Instead, by using the definition of the reduced composite gradient (48), the update rule for the augmented estimate sequence vertices (54) becomes

$$\begin{aligned}\boldsymbol{v}_{k+1} &= \frac{1}{\gamma_{k+1}}\left(\gamma_k\boldsymbol{v}_k - a_{k+1}(L'_{k+1}(\boldsymbol{y}_{k+1} - \boldsymbol{x}_{k+1}) - \mu\boldsymbol{y}_{k+1})\right) \\ &= \frac{1}{\gamma_{k+1}}\left(\gamma_k\boldsymbol{v}_k + a_{k+1}(L_{k+1} + \mu_\Psi)\boldsymbol{x}_{k+1} - a_{k+1}(L_{k+1} - \mu_f)\boldsymbol{y}_{k+1}\right). \end{aligned} \tag{64}$$

Finally, we select the same Armijo-type [1] line-search strategy $\mathcal{S}_A$ as AMGS [13], with parameters $r_u > 1$ and $0 < r_d < 1$ as the increase and, respectively, decrease rates of the Lipschitz constant estimate $\hat{L}_{k+1}$.

In summary, we have established the values of the initial parameters ($A_0 = 0$, $\gamma_0 = 1$, and $\boldsymbol{v}_0 = \boldsymbol{x}_0$), the upper bounds (28), the relaxed supporting generalized parabola lower bounds (39), the line-search strategy $\mathcal{S}_A$, and the expressions of functions $\mathcal{F}_a$ in (63) and $\mathcal{F}_y$ in (61). Based on Algorithm 1, we can now write down ACGM as listed in Algorithm 2.

**Algorithm 2** ACGM in estimate sequence form

**function** ACGM($\boldsymbol{x_0}, L_0, \mu_f, \mu_\Psi, K$)
  $\boldsymbol{v}_0 = \boldsymbol{x}_0,\ \mu = \mu_f + \mu_\Psi,\ A_0 = 0,\ \gamma_0 = 1$ ▷ Initialization
  **for** $k = 0, \ldots, K-1$ **do** ▷ Main loop
    $\hat{L}_{k+1} := r_d L_k$
    **loop**
      $\hat{a}_{k+1} := \frac{1}{2(\hat{L}_{k+1} - \mu_f)}\left(\gamma_k + A_k\mu + \sqrt{(\gamma_k + A_k\mu)^2 + 4(\hat{L}_{k+1} - \mu_f)A_k\gamma_k}\right)$
      $\hat{A}_{k+1} := A_k + \hat{a}_{k+1}$
      $\hat{\gamma}_{k+1} := \gamma_k + \hat{a}_{k+1}\mu$
      $\hat{\boldsymbol{y}}_{k+1} := \frac{1}{A_k\hat{\gamma}_{k+1} + \hat{a}_{k+1}\gamma_k}(A_k\hat{\gamma}_{k+1}\boldsymbol{x}_k + \hat{a}_{k+1}\gamma_k\boldsymbol{v}_k)$
      $\hat{\boldsymbol{x}}_{k+1} := \operatorname{prox}_{\frac{1}{\hat{L}_{k+1}}\Psi}\left(\hat{\boldsymbol{y}}_{k+1} - \frac{1}{\hat{L}_{k+1}}\nabla f(\hat{\boldsymbol{y}}_{k+1})\right)$
      **if** $f(\hat{\boldsymbol{x}}_{k+1}) \leq Q_{f,\hat{L}_{k+1},\hat{\boldsymbol{y}}_{k+1}}(\hat{\boldsymbol{x}}_{k+1})$ **then**
        Break from loop
      **else**
        $\hat{L}_{k+1} := r_u\hat{L}_{k+1}$
      **end if**
    **end loop**
    $L_{k+1} := \hat{L}_{k+1},\ \boldsymbol{x}_{k+1} := \hat{\boldsymbol{x}}_{k+1},\ A_{k+1} := \hat{A}_{k+1},\ \gamma_{k+1} := \hat{\gamma}_{k+1}$
    $\boldsymbol{v}_{k+1} := \frac{1}{\hat{\gamma}_{k+1}}(\gamma_k\boldsymbol{v}_k + \hat{a}_{k+1}(\hat{L}_{k+1} + \mu_\Psi)\hat{\boldsymbol{x}}_{k+1} - \hat{a}_{k+1}(\hat{L}_{k+1} - \mu_f)\hat{\boldsymbol{y}}_{k+1})$
  **end for**
  **return** $x_K$ ▷ Output
**end function**

The per-iteration computational complexity of Algorithm 2 lies between that of FISTA and ACGM. Specifically, an iteration without backtracks requires 1 WTU (as in FISTA / FISTA-CP) while each backtrack costs 2 WTU (see Table 4 for the derivation of these values). Since normal iterations are the most frequent, the average running time behavior of ACGM approaches that of FISTA / FISTA-CP.

**3.4. Convergence analysis.** The convergence of ACGM is governed by (8), with the guarantee given by $A_k$. The growth rate of $A_k$ is affected by the outcome of the line-search procedure. We formulate a simple lower bound for $A_k$ that deals with worst case search behavior. To simplify notation, we introduce the local inverse condition number

$$q_{k+1} \stackrel{\text{def}}{=} \frac{\mu}{L'_{k+1}} = \frac{\mu}{L_{k+1} + \mu_\Psi}, \quad k \geq 0. \tag{65}$$

Table 4: Resource allocation and runtime behavior of parallel black-box ACGM

| Situation | WTU | PPU 1 | | PPU 2 | | PPU 3 | |
|---|---|---|---|---|---|---|---|
| | | Comp. | Iter. | Comp. | Iter. | Comp. | Iter. |
| No backtracks | T | $\nabla f(\boldsymbol{y}_{k+1})$ | k | $f(\boldsymbol{y}_{k+1})$ | k | $f(\boldsymbol{x}_k)$ | k - 1 |
| | T + 1 | $\nabla f(\boldsymbol{y}_{k+2})$ | k + 1 | $f(\boldsymbol{y}_{k+2})$ | k + 1 | $f(\boldsymbol{x}_{k+1})$ | k |
| | T + 2 | $\nabla f(\boldsymbol{y}_{k+3})$ | k + 2 | $f(\boldsymbol{y}_{k+3})$ | k + 2 | $f(\boldsymbol{x}_{k+2})$ | k + 1 |
| Single backtrack | T | $\nabla f(\boldsymbol{y}_{k+1})$ | k | $f(\boldsymbol{y}_{k+1})$ | k | $f(\boldsymbol{x}_k)$ | k - 1 |
| | T + 1 | $\nabla f(\boldsymbol{y}_{k+2})$ | k + 1 | $f(\boldsymbol{y}_{k+2})$ | k + 1 | $f(\boldsymbol{x}_{k+1})$ | k |
| | T + 2 | $\nabla f(\boldsymbol{y}_{k+1})$ | k | $f(\boldsymbol{y}_{k+1})$ | k | Idle | |
| | T + 3 | $\nabla f(\boldsymbol{y}_{k+2})$ | k + 1 | $f(\boldsymbol{y}_{k+2})$ | k + 1 | $f(\boldsymbol{x}_{k+1})$ | k |
| | T + 4 | Idle | | Idle | | $f(\boldsymbol{x}_{k+2})$ | k + 1 |

If $L_{k+1} \geq L_f$, then the descent condition for $f$ (27) holds regardless of the algorithmic state, implying the backtracking search will guarantee that

$$L_{k+1} \leq L_u \stackrel{\text{def}}{=} \max\{r_u L_f, r_d L_0\}, \quad k \geq 0. \tag{66}$$

Let the worst case local inverse condition number be defined as

$$q_u \stackrel{\text{def}}{=} \frac{\mu}{L_u + \mu_\Psi} \leq q_{k+1}, \quad k \geq 0. \tag{67}$$

THEOREM 4. *The convergence guarantee $A_k$ for ACGM in the non-strongly convex case ($\mu = 0$) is lower bounded by*

$$A_k \geq \frac{(k+1)^2}{4L_u}, \quad k \geq 1, \tag{68}$$

*and in the strongly convex case ($\mu > 0$) by*

$$A_k \geq \frac{1}{L_u - \mu_f}(1 - \sqrt{q_u})^{-(k-1)}, \quad k \geq 1. \tag{69}$$

*Proof.* In the non-strongly convex case, we have

$$\begin{aligned} A_{k+1} = A_k + a_{k+1} &\stackrel{(60)}{\geq} A_k + \frac{1 + \sqrt{1 + 4L_{k+1}A_k}}{2L_{k+1}} \\ &\stackrel{(66)}{\geq} A_k + \frac{1}{2L_u} + \sqrt{\frac{1}{4L_u^2} + \frac{A_k}{L_u}}. \end{aligned} \tag{70}$$

We prove by induction that (68) holds for all $k \geq 1$. First, for $k = 1$, (68) is valid since

$$A_1 = \frac{1}{L_1} \geq \frac{(1+1)^2}{4L_u}. \tag{71}$$

Next, we assume that (68) is valid for $k$, and show that it holds for $k + 1$. From (70) and (68), we have

$$\begin{aligned} A_{k+1} &\geq \frac{1}{4L_u}\left((k+1)^2 + 2 + 2\sqrt{1 + (k+1)^2}\right) \\ &\geq \frac{(k+2)^2}{4L_u}. \end{aligned} \tag{72}$$

In the strongly convex case, the curvature of the estimate function can be expressed in absolute terms as

$$\gamma_k = \gamma_0 + \left(\sum_{i=1}^{k} a_i\right)\mu = \gamma_0 + (A_k - A_0)\mu = 1 + A_k\mu, \quad k \geq 0, \tag{73}$$

which trivially implies that $\gamma_k > A_k\mu$. Hence, we have

$$\frac{a_{k+1}^2}{A_{k+1}^2} \overset{(62)}{=} \frac{\gamma_{k+1}}{(L_{k+1} + \mu_\Psi)A_{k+1}} > \frac{\mu}{L_{k+1} + \mu_\Psi} = q_{k+1} \geq q_u, \quad k \geq 0, \tag{74}$$

which leads to

$$\frac{A_{k+1}}{A_k} > \frac{1}{1 - \sqrt{q_u}}, \quad k \geq 1. \tag{75}$$

Using $A_1 = \frac{1}{L_1 - \mu_f} \geq \frac{1}{L_u - \mu_f}$, the strongly convex lower bound (69) follows by induction. ∎

**3.5. ACGM in extrapolated form.** An interesting property of FGM is that, at every iteration $k$, the next point where the gradient is queried $\boldsymbol{y}_{k+2}$ can be expressed in terms of the previous two iterates $\boldsymbol{x}_{k+1}$ and $\boldsymbol{x}_k$ by extrapolation, namely

$$\boldsymbol{y}_{k+2} = \boldsymbol{x}_{k+1} + \beta_{k+1}(\boldsymbol{x}_{k+1} - \boldsymbol{x}_k), \quad k \geq 0, \tag{76}$$

where $\beta_{k+1}$ is the auxiliary point extrapolation factor.

We demonstrate that this applies to ACGM (Algorithm 2) as well, with the difference that $\beta_{k+1}$ can only be computed during iteration $k+1$ due to uncertainty in the outcome of line-search. To be able to use extrapolation (76) in the first iteration, we define $\boldsymbol{x}_{-1} \overset{\text{def}}{=} \boldsymbol{x}_0$. First, we show the following property of ACGM, which carries over from FGM.

LEMMA 5. *The estimate function vertices can be obtained from successive iterates through extrapolation as*

$$\boldsymbol{v}_{k+1} = \boldsymbol{x}_k + \frac{A_{k+1}}{a_{k+1}}(\boldsymbol{x}_{k+1} - \boldsymbol{x}_k), \quad k \geq 0. \tag{77}$$

*Proof.* By combining (59) with (64), we get

$$\begin{aligned}
\boldsymbol{v}_{k+1} &= \frac{\gamma_k}{\gamma_{k+1}} \frac{(a_{k+1}\gamma_k + A_k\gamma_{k+1})\boldsymbol{y}_{k+1} - A_k\gamma_{k+1}\boldsymbol{x}_k}{a_{k+1}\gamma_k} + \\
&\quad + \frac{a_{k+1}(L_{k+1} + \mu_\Psi)}{\gamma_{k+1}}\boldsymbol{x}_{k+1} - \frac{a_{k+1}(L_{k+1} - \mu_f)}{\gamma_{k+1}}\boldsymbol{y}_{k+1} \\
&\overset{(62)}{=} \frac{a_{k+1}\gamma_k + A_k\gamma_{k+1} - A_{k+1}\gamma_{k+1} - a_{k+1}^2\mu}{a_{k+1}\gamma_{k+1}}\boldsymbol{y}_{k+1} + \frac{A_{k+1}}{a_{k+1}}\boldsymbol{x}_{k+1} - \frac{A_k}{a_{k+1}}\boldsymbol{x}_k \\
&\overset{(53)}{=} \boldsymbol{x}_k + \frac{A_{k+1}}{a_{k+1}}(\boldsymbol{x}_{k+1} - \boldsymbol{x}_k), \quad k \geq 0. \qquad (78)
\end{aligned}$$

∎

Now, from (59) and Lemma 5, we obtain that the extrapolation expression (76) does indeed hold for ACGM, namely

$$\begin{aligned}
\boldsymbol{y}_{k+1} &= \frac{1}{A_k\gamma_{k+1} + a_{k+1}\gamma_k}\left(A_k\gamma_{k+1}\boldsymbol{x}_k + a_{k+1}\gamma_k\left(\boldsymbol{x}_{k-1} + \frac{A_k}{a_k}(\boldsymbol{x}_k - \boldsymbol{x}_{k-1})\right)\right) \\
&= \boldsymbol{x}_k + \beta_k(\boldsymbol{x}_k - \boldsymbol{x}_{k-1}), \quad k \geq 0, \qquad (79)
\end{aligned}$$

where the auxiliary point extrapolation factor $\beta_k$ is given by

$$\beta_k = \frac{a_{k+1}\gamma_k\left(\frac{A_k}{a_k} - 1\right)}{A_k\gamma_{k+1} + a_{k+1}\gamma_k}, \quad k \geq 0. \tag{80}$$

To bring ACGM to a form in which it can be easily compared with FISTA / FISTA-CP, we denote the vertex extrapolation factor in Lemma 5 as

$$t_k \stackrel{\text{def}}{=} \begin{cases} \frac{A_k}{a_k}, & k \geq 1, \\ 0, & k = 0. \end{cases} \tag{81}$$

The accumulated weights and the curvature ratio $\gamma_k/\gamma_{k+1}$ can be written in terms of $t_k$ as

$$\begin{aligned} A_{k+1} &\stackrel{(62)}{=} \frac{A_{k+1}^2\gamma_{k+1}}{(L_{k+1}+\mu_\Psi)a_{k+1}^2} \stackrel{(81)}{=} \frac{\gamma_{k+1}t_{k+1}^2}{L_{k+1}+\mu_\Psi}, \quad k \geq 0, \\ A_0 &= 0 \stackrel{(81)}{=} \frac{\gamma_0 t_0^2}{L_0+\mu_\Psi}, \end{aligned} \tag{82}$$

$$\frac{\gamma_k}{\gamma_{k+1}} = 1 - \frac{A_{k+1}a_{k+1}\mu}{A_{k+1}\gamma_{k+1}} \stackrel{(62)}{=} 1 - \frac{A_{k+1}a_{k+1}\mu}{(L_{k+1}+\mu_\Psi)a_{k+1}^2} \stackrel{(81)}{=} 1 - q_{k+1}t_{k+1}, \quad k \geq 0. \tag{83}$$

Expressions (82) and (83) facilitate the derivation of a recursion rule for $t_k$ that does not depend on either $a_k$ or $A_k$, as follows:

$$(L_{k+1}+\mu_\Psi)A_{k+1} - (L_{k+1}+\mu_\Psi)a_{k+1} - \frac{L_{k+1}+\mu_\Psi}{L_k+\mu_\Psi}(L_k+\mu_\Psi)A_k = 0 \quad \stackrel{(82)}{\Leftrightarrow} \tag{84}$$

$$\gamma_{k+1}t_{k+1}^2 - \gamma_{k+1}t_{k+1} - \frac{L_{k+1}+\mu_\Psi}{L_k+\mu_\Psi}\gamma_k t_k^2 = 0 \quad \stackrel{(83)}{\Leftrightarrow} \tag{85}$$

$$t_{k+1}^2 + t_{k+1}(q_k t_k^2 - 1) - \frac{L_{k+1}+\mu_\Psi}{L_k+\mu_\Psi}t_k^2 = 0, \quad k \geq 0, \quad \mu \geq 0. \tag{86}$$

Lastly, we write down the auxiliary point extrapolation factor $\beta_k$ as

$$\begin{aligned} \beta_k &= \frac{a_{k+1}\gamma_k\left(\frac{A_k}{a_k} - 1\right)}{A_k\gamma_{k+1} + a_{k+1}\gamma_k} \stackrel{(81)}{=} \frac{t_k - 1}{t_{k+1}} \frac{A_{k+1}\gamma_k}{A_k\gamma_{k+1} + a_{k+1}\gamma_k} \\ &= \frac{t_k - 1}{t_{k+1}} \frac{\frac{\gamma_k}{\gamma_{k+1}}}{\frac{A_{k+1}\gamma_{k+1} - a_{k+1}(\gamma_{k+1}-\gamma_k)}{A_{k+1}\gamma_{k+1}}} \stackrel{(53)}{=} \frac{t_k - 1}{t_{k+1}} \frac{\frac{\gamma_k}{\gamma_{k+1}}}{1 - \frac{\mu a_{k+1}^2}{A_{k+1}\gamma_{k+1}}} \\ &\stackrel{(83)}{=} \frac{t_k - 1}{t_{k+1}} \frac{1 - q_{k+1}t_{k+1}}{1 - q_{k+1}}, \quad k \geq 0. \end{aligned} \tag{87}$$

Now, from (86) and (87), we can formulate ACGM based on extrapolation, as presented in Algorithm 3. Algorithms 2 and 3 differ in form but are theoretically guaranteed to produce identical iterates.

## 4. Numerical results.

In this section we test ACGM against the state-of-the-art methods on two problems representative of their respective classes. In subsection 4.1 we show a typical non-strongly convex application whereas in subsection 4.2 we focus on a strongly convex problem.

**Algorithm 3** ACGM in extrapolated form

**function** ACGM($\boldsymbol{x_0}, L_0, \mu_f, \mu_\Psi, K$)
  $\boldsymbol{x}_{-1} = \boldsymbol{x}_0,\ \mu = \mu_f + \mu_\Psi,\ t_0 = 0,\ q_0 = \frac{\mu}{L_0+\mu_\Psi}$ ▷ Initialization
  **for** k = 0,...,K-1 **do** ▷ Main loop
    $\hat{L}_{k+1} := r_d L_k$
    **loop**
      $\hat{q}_{k+1} := \frac{\mu}{\hat{L}_{k+1}+\mu_\Psi}$
      $\hat{t}_{k+1} := \frac{1}{2}\left(1 - q_k t_k^2 + \sqrt{(1-q_k t_k^2)^2 + 4\frac{\hat{L}_{k+1}+\mu_\Psi}{L_k+\mu_\Psi} t_k^2}\right)$
      $\hat{\boldsymbol{y}}_{k+1} := \boldsymbol{x}_k + \frac{t_k-1}{\hat{t}_{k+1}} \frac{1-\hat{q}_{k+1}\hat{t}_{k+1}}{1-\hat{q}_{k+1}} (\boldsymbol{x}_k - \boldsymbol{x}_{k-1})$
      $\hat{\boldsymbol{x}}_{k+1} := \text{prox}_{\frac{1}{\hat{L}_{k+1}}\Psi}\left(\hat{\boldsymbol{y}}_{k+1} - \frac{1}{\hat{L}_{k+1}}\nabla f(\hat{\boldsymbol{y}}_{k+1})\right)$
      **if** $f(\hat{\boldsymbol{x}}_{k+1}) \le Q_{f,\hat{L}_{k+1},\hat{\boldsymbol{y}}_{k+1}}(\hat{\boldsymbol{x}}_{k+1})$ **then**
        Break from loop
      **else**
        $\hat{L}_{k+1} := r_u \hat{L}_{k+1}$
      **end if**
    **end loop**
    $\boldsymbol{x}_{k+1} := \hat{\boldsymbol{x}}_{k+1},\ L_{k+1} := \hat{L}_{k+1},\ q_{k+1} := \hat{q}_{k+1},\ t_{k+1} := \hat{t}_{k+1}$
  **end for**
  **return** $x_K$ ▷ Output
**end function**

**4.1. L1-regularized image deblurring.** We choose to test the capabilities of ACGM (Algorithm 3) on the very problem FISTA was introduced to solve, namely the $l_1$ regularized deblurring of images. We have also noticed that several monographs in the field (e.g. [6, 17]) do not include AMGS in their benchmarks. For completeness, we compare Algorithm 3 against both FISTA with backtracking line-search and AMGS. For ease and accuracy of benchmarking, we have adopted the experimental setup from Section 5.1 in [2]. Here, the composite objective function is given by

$$f(\boldsymbol{x}) = \|\boldsymbol{A}\boldsymbol{x} - \boldsymbol{b}\|_2^2, \qquad \Psi(\boldsymbol{x}) = \lambda\|\boldsymbol{x}\|_1, \tag{88}$$

where $\boldsymbol{A} = \boldsymbol{R}\boldsymbol{W}$. Linear operator $\boldsymbol{R}$ is Gaussian blur with standard deviation 4.0 and $9 \times 9$ pixel kernel, applied using reflexive boundary conditions [7]. Linear operator $\boldsymbol{W}$ is the inverse three-stage Haar wavelet transform. Variable $\boldsymbol{x} \in \mathbb{R}^{n_1 \times n_2}$ is a digital image of dimensions $n_1 = n_2 = 256$. The blurred image $\boldsymbol{b}$ is obtained by applying $\boldsymbol{R}$ to the cameraman test image [2] with pixel values scaled to the $[0, 1]$ range, followed by the addition of Gaussian noise (zero-mean, standard deviation $10^{-3}$). The constant $L_f$ can be computed as the maximum eigenvalue of a symmetric Toeplitz-plus-Hankel matrix (more details in [7]), which yields a value of $L_f = 2.0$. The problem is non-strongly convex with $\mu = \mu_f = \mu_\Psi = 0$. The regularization parameter $\lambda$ is set to $2 \cdot 10^{-5}$ to adjust for the noise level of $\boldsymbol{b}$.

The starting point $\boldsymbol{x}_0$ was set to $\boldsymbol{W}^{-1}\boldsymbol{b}$ for all algorithms. AMGS and FISTA were run using $r_u^{\text{AMGS}} = r_u^{\text{FISTA}} = 2.0$ and $r_d^{\text{AMGS}} = 0.9$ as these values were suggested in [3] to "provide good performance in many applications". Assuming that most of time the Lipschitz constant estimates hover around a fixed value, we have for AMGS that a backtrack occurs every $-(\log r_u^{\text{AMGS}})/(\log r_d^{\text{AMGS}})$ iterations. Since the cost ratio between a backtrack and non-backtrack iteration for the proposed method,

ACGM, is double that of AMGS, to ensure that the line-search procedures of both methods have comparable computational overhead, we have chosen $r_d^{\text{ACGM}} = r_d^{\text{AMGS}}$ and $r_d^{\text{ACGM}} = \sqrt{r_d^{\text{AMGS}}}$.

The line-search procedure of FISTA differs from that of ACGM in that the Lipschitz constant estimates can only increase while the algorithm is running. To showcase the robustness of the line-search employed by ACGM (and AMGS), we have considered two scenarios: a pathologically overestimated initial guess $L_0 = 10L_f$ and a normally underestimated $L_0 = 0.3L_f$ (Figure 1). The convergence rate is measured as the difference between objective function values and an optimal value *estimate* $F(\hat{\boldsymbol{x}}^*)$, where $\hat{\boldsymbol{x}}^*$ is the iterate obtained after running fixed step size FISTA with the correct Lipschitz constant parameter for 10000 iterations.

When indexed in iterations (Figures 1a and 1b), ACGM converges roughly as fast as AMGS. ACGM takes the lead after 500 iterations, due to the stringency of AMGS's "damped relaxation condition" [13] compared to the ACGM's descent condition. However, when indexed in WTU, ACGM clearly surpasses AMGS (Figures 1c and 1d), the main reason being ACGM's low per-iteration complexity. FISTA lags behind in the overestimated case, regardless of the convergence measure (Figures 1a and 1c), and it is also slightly slower in the underestimated case (Figure 1d). The disadvantages of FISTA's line-search are evidenced by Figures 1e and 1f. In both cases, FISTA produces, on average, a higher Lipschitz estimate than ACGM.

**4.2. Total variation based image denoising.** The treatment of first-order algorithms for strongly convex objectives is relatively rare in the literature. The recent monograph [6] provides an extensive review of such algorithms (in example 4.14) on an image restoration problem, specifically the problem of minimizing the dual of the Huber-Rudin-Osher-Fatemi (Huber-ROF) [6, 18], a variant of the total variation based image denoising problem. Again, we note there the lack of comparison against AMGS and choose, for ease of benchmarking, this same application to test Algorithm 3. We offer a brief description of the problem and forward the reader to [6] for details.

The objective function is given by

$$
\begin{aligned}
f(\boldsymbol{p}) &= \frac{1}{2}\|\boldsymbol{D}^*\boldsymbol{p} - \boldsymbol{b}\|_2^2, \\
\Psi(\boldsymbol{p}) &= \begin{cases} \frac{\epsilon}{2\lambda}\|\boldsymbol{p}\|_2^2 & \text{if } \|\boldsymbol{p}_{i,j}\|_2^2 \le \lambda \ \text{ for all } i = 1,\dots,n_1,\ j = 1,\dots,n_2, \\ \infty & \text{otherwise.} \end{cases}
\end{aligned} \tag{89}
$$

Here $\boldsymbol{p} \in \mathbb{R}^{n_1\times n_2\times 2}$ and $n_1 = n_2 = 256$. The discrete gradient operator is defined as

$$
\begin{aligned}
&\boldsymbol{D} : \mathbb{R}^{n_1\times n_2} \to \mathbb{R}^{n_1\times n_2\times 2}, \\
&(\boldsymbol{D}\boldsymbol{u})_{i,j} = \left(\boldsymbol{u}_{\min\{i+1,n_1\},j} - \boldsymbol{u}_{i,j}, \boldsymbol{u}_{i,\min\{j+1,n_2\}} - \boldsymbol{u}_{i,j}\right),
\end{aligned} \tag{90}
$$

with $\boldsymbol{D}^*$ as its adjoint. The $n_1\times n_2$ real-valued digital image $\boldsymbol{b}$ is obtained by applying Gaussian noise (zero-mean, standard deviation 0.1) to the cameraman test image [2] with pixel values scaled to the $[0,1]$ range. Function $f$ has Lipschitz gradient with $L_f \le 8$ (see [4]) and $\mu_f = 0$. In function $\Psi$, $\epsilon = 0.001$ is a primal smoothing parameter (see also [12]) and $\lambda = 0.1$ is a regularization parameter. Hence $\mu = \mu_\Psi = 0.01$.

We benchmark ACGM against methods equipped with a line-search procedure, such as FISTA and AMGS, and with methods that rely on $L_f$ being known in advance, namely a proximal variant of Nesterov's Scheme III [6, 11] and FISTA-CP. The starting point for all algorithms was set to $\boldsymbol{x}_0 = \boldsymbol{D}\boldsymbol{b}$. For the same reasons outlined

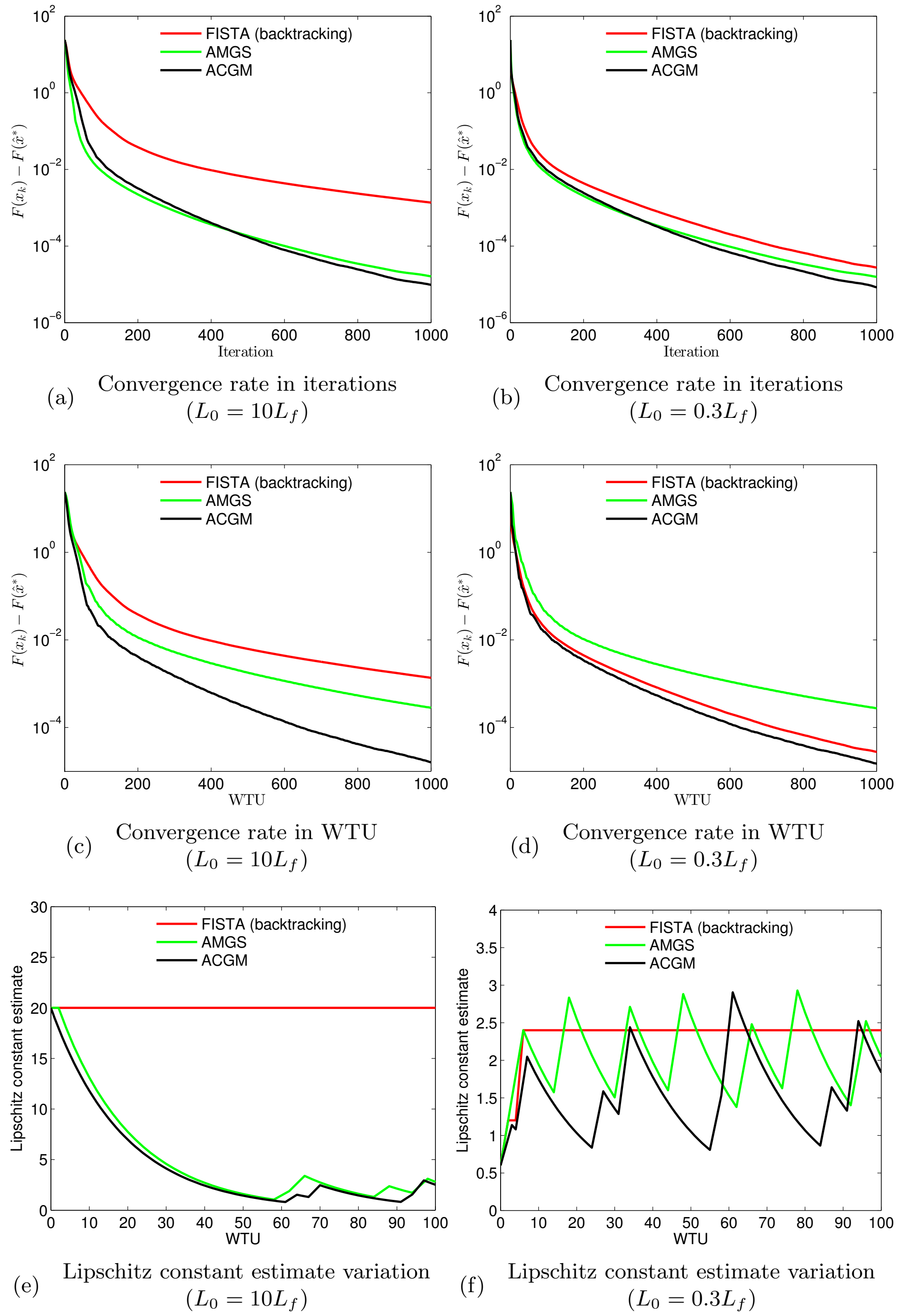

(a) Convergence rate in iterations ($L_0 = 10L_f$)

(b) Convergence rate in iterations ($L_0 = 0.3L_f$)

(c) Convergence rate in WTU ($L_0 = 10L_f$)

(d) Convergence rate in WTU ($L_0 = 0.3L_f$)

(e) Lipschitz constant estimate variation ($L_0 = 10L_f$)

(f) Lipschitz constant estimate variation ($L_0 = 0.3L_f$)

Fig. 1: Convergence results on the L1-regularized image deblurring problem ($\mu = 0$)

in subsection 4.1, we have chosen $r_u^{\mathrm{ACGM}} = r_u^{\mathrm{AMGS}} = r_u^{\mathrm{FISTA}} = 2.0$, $r_d^{\mathrm{AMGS}} = 0.9$, and $r_d^{\mathrm{ACGM}} = \sqrt{r_d^{\mathrm{AMGS}}}$.

Our system had limited floating point precision and the descent condition stopping criterion in oracle form (Algorithm 3, line 10) was be prone to numerical errors. It often produced values of $L_{k+1}$ beyond $L_u$, or failed to terminate the search loop altogether. By exploiting the quadratic structure of $f$, we have replaced the descent condition stopping criterion with an equivalent and more numerically stable condition, written as

$$\|D^*(\boldsymbol{x}_{k+1} - \boldsymbol{y}_{k+1})\|_2^2 \leq L_{k+1}\|\boldsymbol{x}_{k+1} - \boldsymbol{y}_{k+1}\|_2^2, \quad k \geq 0. \tag{91}$$

We have chosen not to add the proximal heavy ball method (PHBM) to our benchmark (as in [6]) because it lacks formal convergence guarantees and its surprising performance on this quadratically constrained quadratic program (QCQP) does not generalize to the entire composite objective problem class. However, we have employed PHBM to compute the optimal point estimate $\hat{\boldsymbol{x}}^*$ as the iterate obtained after running PHBM with $L_f = 8$ for 800 iterations. Methods equipped with a line-search procedure incur a search overhead whereas the other methods do not. For fair comparison, we have tested the collection of methods in the accurate $L_f = 8$ case as well as the lossy $L_f = 20$ case (Figure 2).

When indexing in iterations, AMGS converges the fastest (Figures 2a and 2b) although it is inferior to ACGM and FISTA-CP in terms of WTU usage (Figures 2c and 2d). FISTA with backtracking is unable to reach even a linear convergence rate. It leads only during the first 100 iterations when an accurate value of $L_f$ is supplied (Figure 2c) but lags behind afterwards, without triggering any function scheme adaptive restart condition [16]. Overall, we deem it unsuitable for strongly convex objectives (Figures 2c and 2d). Nesterov's scheme III displays a comparable asymptotic rate to ACGM when $L_f$ is accurate (Figures 2a and 2c) but slower convergence altogether (Figures 2a to 2d). Even though it incurs a search overhead, ACGM narrowly outperformes FISTA-CP in the accurate case (Figure 2c) by exploiting the fact that the local curvature is often below $L_f$ (Figure 2e). In the inaccurate case (Figure 2d), ACGM leads all other algorithms by a large margin due to the efficiency of the line-search procedure (Figure 2f), which easily compensates for the increased overhead.

**5. Discussion.** The proposed method, ACGM, when formulated using extrapolation, encompasses several existing optimization schemes. Specifically, Algorithm 3 *without the line-search procedure*, i.e., with $L_k = L_f$ for all $k \geq 0$, produces the same iterates as FISTA-CP with the theoretically optimal step size $\tau^{\mathrm{FISTA-CP}} = \frac{1}{L_f}$. In the non-strongly convex case, ACGM without the line-search reduces to constant step size FISTA. Also for $\mu = 0$, ACGM with line-search constitutes a simplified and more intuitive alternative to a recently introduced (without derivation) line-search extension of FISTA [19].

However, ACGM is more than an umbrella method and actually surpasses constituent FISTA and FISTA-CP as well as the adaptive AMGS and the efficient FGM. Indeed, FISTA suffers from two drawbacks: the parameter $t_k^{\mathrm{FISTA}}$ update is oblivious to the change in local curvature and the Lipschitz constant estimates cannot decrease. Hence, if the initial Lipschitz estimate is erroneously large, FISTA will slow down considerably (see also subsection 4.1). We formally express the advantages of ACGM's line-search over that of FISTA in the following proposition.

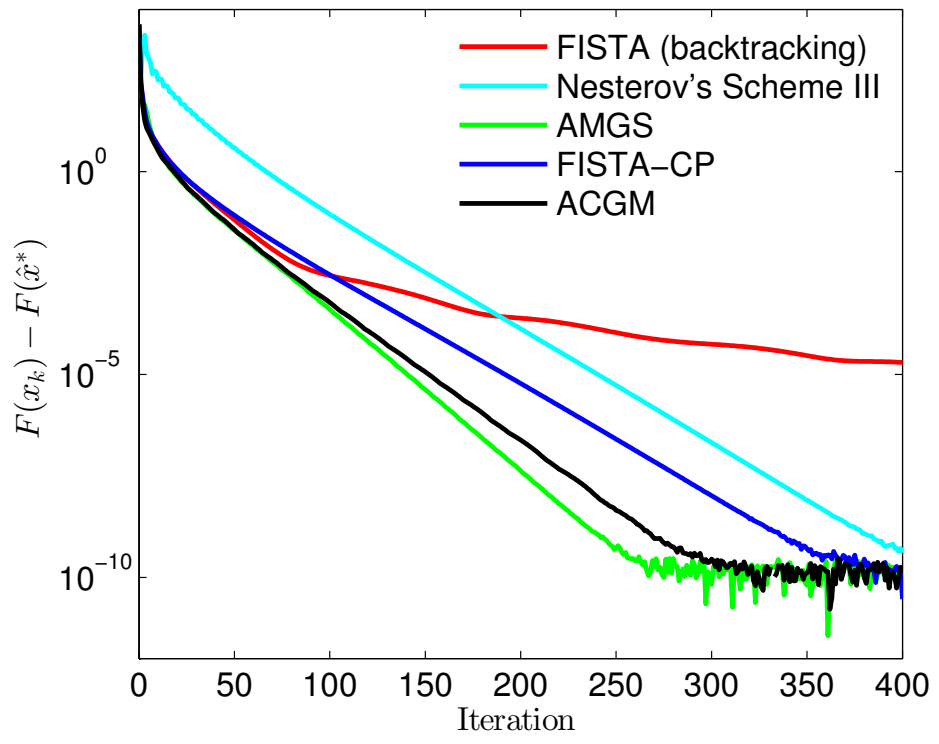

(a) Convergence rate in iterations ($L_f = 8$)

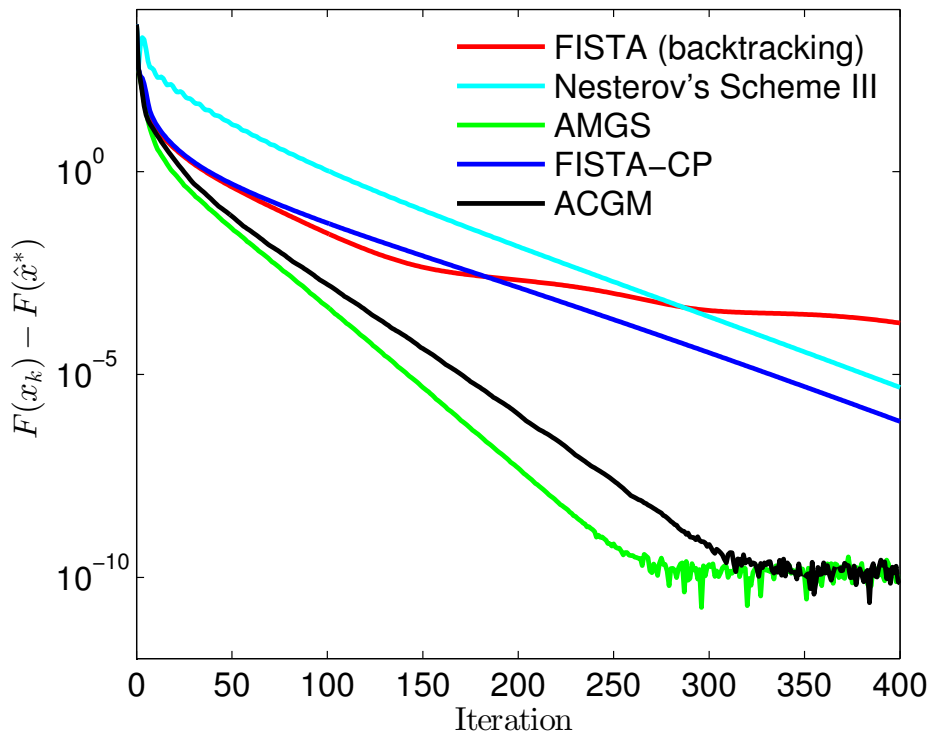

(b) Convergence rate in iterations ($L_f = 20$)

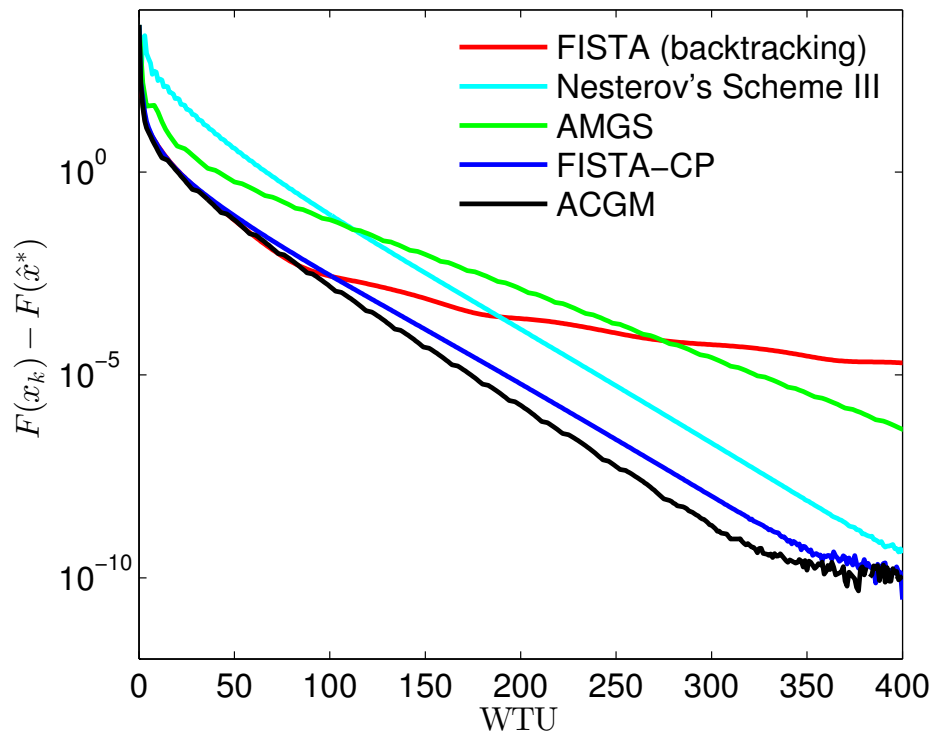

(c) Convergence rate in WTU ($L_f = 8$)

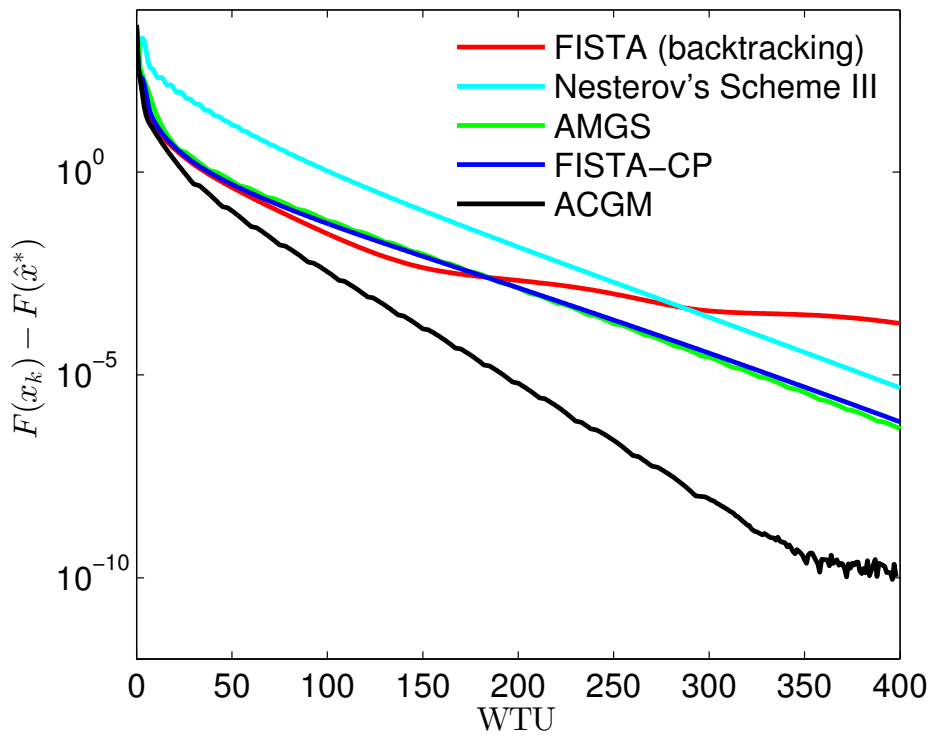

(d) Convergence rate in WTU ($L_f = 20$)

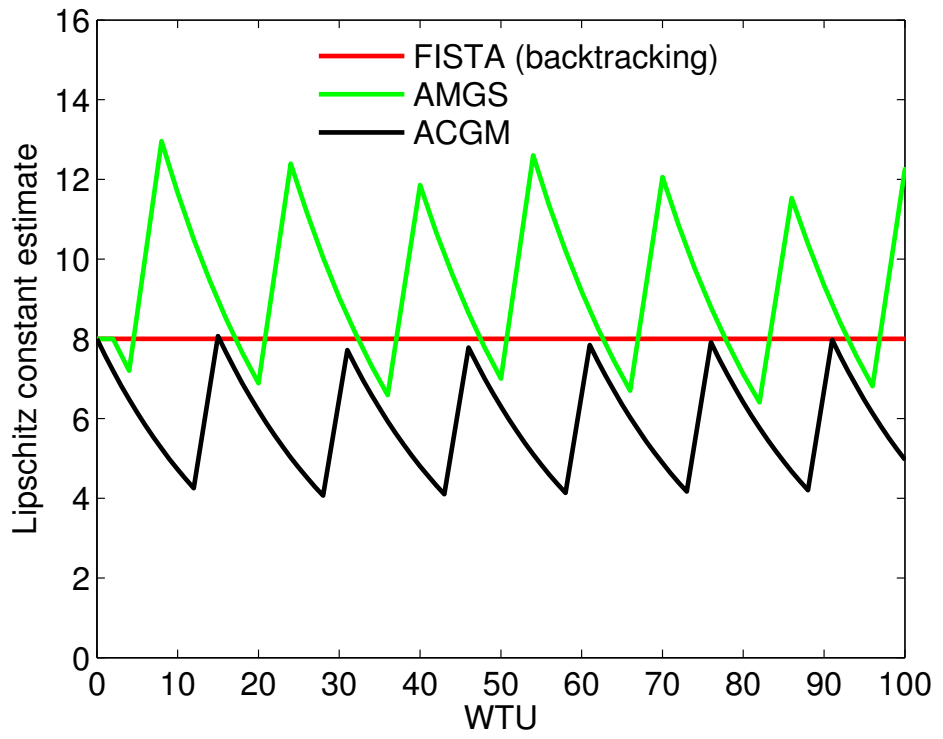

(e) Lipschitz constant estimate variation ($L_f = 8$)

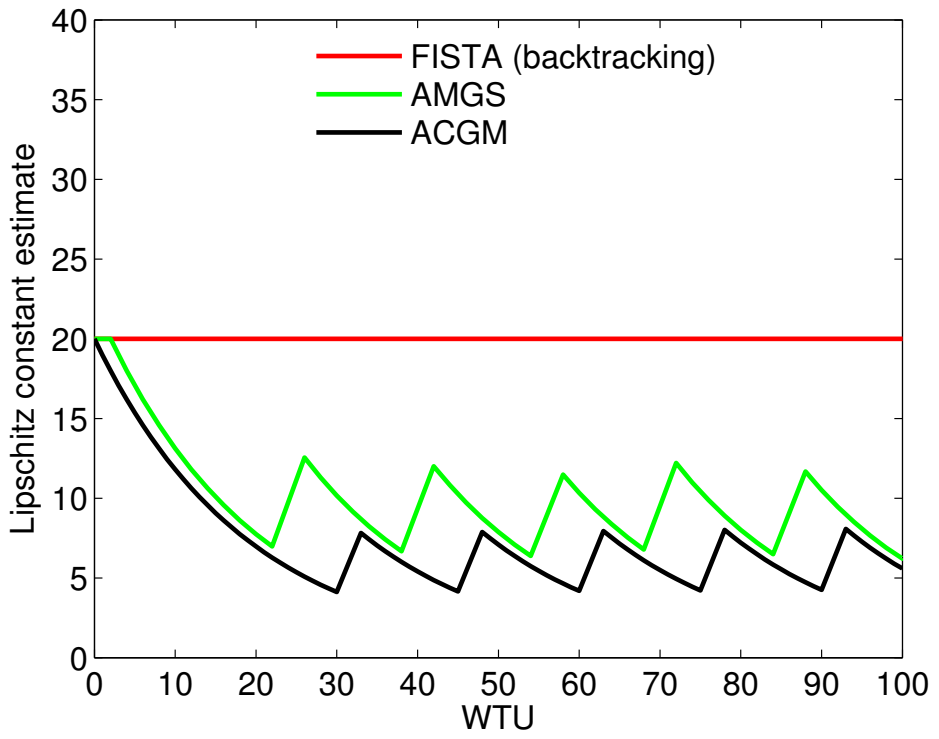

(f) Lipschitz constant estimate variation ($L_f = 20$)

Fig. 2: Convergence results on the dual Huber-ROF problem ($\mu = 0.01$)

PROPOSITION 6. *In the non-strongly convex case* ($\mu = 0$), *under identical local curvature conditions, when* $r_u^{\text{ACGM}} = r_u^{\text{FISTA}}$, *ACGM has superior theoretical convergence guarantees to FISTA, namely*

$$A_k^{\text{ACGM}} \geq A_k^{\text{FISTA}}, \quad k \geq 0. \tag{92}$$

*Proof.* See Appendix C. □

The ability to dynamically and frequently adjust to the local Lipschitz constant gives ACGM an advantage over FISTA-CP, even when an accurate estimate of the Lipschitz constant is available beforehand. ACGM has the provable convergence rate of FGM (Theorem 4). Hence, by expanding the arguments presented in Appendix B, ACGM is theoretically guaranteed to outperform AMGS as well.

The per-iteration complexity of ACGM, both in the non-strongly and strongly convex cases ($\mu \geq 0$), lies well below that of AMGS. Considering that backtracks rarely occur, it approaches that of FISTA (see Table 5). This is close to smallest possible value of 1 WTU per iteration. All of the above theoretical performance guarantees are corroborated by the simulation results in section 4.

Table 5: Iteration cost in WTU of line-search methods AMGS, FISTA, and ACGM

| Iteration phase | AMGS | FISTA | ACGM |
|---|---|---|---|
| Iteration without backtrack | 2 | 1 | 1 |
| Each backtrack | 2 | 1 | 2 |

Thus, this is the first time, as far as we are aware, that a method has been shown to be superior, from theoretical as well as simulation results, to both AMGS and FISTA / FISTA-CP. ACGM is also the most robust first-order method among the state-of-the-art for its problem class, in the sense that it is the only method that is able to deliver the speed of FGM while being as widely applicable as AMGS (see Table 6 for detailed feature comparison).

Table 6: Features of black-box first-order methods

| Feature | Prox. point | FGM | AMGS | FISTA | FISTA-CP | ACGM |
|---|---|---|---|---|---|---|
| Non-differentiable objective | yes | no | yes | yes | yes | yes |
| Line-search | no | no | yes | partial | no | yes |
| $\mathcal{O}(\frac{1}{k^2})$ rate for $\mu = 0$ | no | yes | yes | yes | yes | yes |
| Linear rate for $\mu > 0$ | yes | yes | yes | no | yes | yes |
| $\mathcal{O}((1-\sqrt{q})^k)$ rate for $\mu > 0$ | no | yes | no | no | yes | yes |

All these capabilities stem from the augmented estimate sequence framework. The results presented in this work confirm that the estimate sequence is an effective tool for designing fast algorithms and suggest that this concept can be extended to problems outside its original scope.

**Appendix A. FGM estimate sequences.**

When constructing FGM for differentiable (and hence finite) objectives, Nesterov expresses the convergence guarantee (8) in [11] as

$$A_k(F(\boldsymbol{x}_k) - F(\boldsymbol{x}^*)) \le A_0(F(\boldsymbol{x}_0) - F(\boldsymbol{x}^*)) + \frac{\gamma_0}{2}\|\boldsymbol{x}_0 - \boldsymbol{x}^*\|_2^2, \quad x^* \in X^*, \quad k \ge 0.$$

The highest allowable upper bound on weighted function values $A_k F(\boldsymbol{x}_k)$ becomes

$$H_k \overset{\text{def}}{=} (A_k - A_0)F(\boldsymbol{x}^*) + A_0 F(\boldsymbol{x}_0) + \frac{\gamma_0}{2}\|\boldsymbol{x}_0 - \boldsymbol{x}^*\|_2^2.$$

We interpret the estimate sequence definition in [11] as

$$\psi_k(\boldsymbol{x}) \overset{\text{def}}{=} (A_k - A_0)W_k(\boldsymbol{x}) + A_0 F(\boldsymbol{x}_0) + \frac{\gamma_0}{2}\|\boldsymbol{x} - \boldsymbol{x}_0\|_2^2, \quad k \ge 0,$$

where $W_k$ are global convex lower bounds on $F$ and $W_0$ is unspecified. The convergence guarantee follows in the same way as in subsection 2.1, that is,

$$A_k F(\boldsymbol{x}_k) \le \psi_k^* \le \psi_k(\boldsymbol{x}^*) \le H_k, \quad k \ge 0.$$

It is clear that if $A_0$ is constrained to be positive, it can be set to any positive number without changing the convergence guarantees, provided that $\gamma_0$ along with $A_k$ for $k \ge 1$ are scaled accordingly. We choose without loss of generality $A_0 = 1$. As argued in [11], we have $\gamma_0 \ge \mu$ and $\gamma_0 > 0$ without the restriction $\gamma_0 \le 1$.

**Appendix B. Proving that FGM converges faster than AMGS.**

To be able to compare the provable convergence rates of FGM and AMGS, we consider the largest problem class to which both algorithms are applicable. Let $\Psi(\boldsymbol{x}) = 0$ and $L_f$ known in advance. For ease of analysis, we study AMGS without line-search (and the corresponding search overhead). Note that FGM and AMGS are not parallelizable in black-box form. Hence, the comparison carries over to the serial scenario as well.

In the non-strongly convex case, the convergence guarantee for AMGS is given by

$$A_k^{\text{AMGS}} = A_{\frac{i}{2}}^{\text{AMGS}} \ge \frac{i^2}{8L_f},$$

where $i$ gives the number of WTU required by the first $k$ iterations. This is asymptotically identical to the rate of FGM for a conservative parameter choice of $\gamma_0 = L_f$ that yields

$$A_k^{\text{FGM}} = A_i^{\text{FGM}} \ge \frac{(i+2)^2}{8L_f}.$$

In the strongly convex case, let $q$ be the inverse condition number of the objective function, $q \overset{\text{def}}{=} \frac{\mu}{L_f}$. We assume $q < 1$ since $q = 1$ means that the optimization problem can be solved exactly, using only one iteration of either AMGS or FGM. When employing AMGS, Nesterov suggests in [13] either to transfer all strong convexity from $f$ to $\Psi$, or to restart the algorithm at regular intervals. Both enhancements have the same effect on the convergence guarantee, which can be expressed as

$$A_k^{\text{AMGS}} = A_{\frac{i}{2}}^{\text{AMGS}} \ge C^{\text{AMGS}} \left(B^{\text{AMGS}}\right)^i,$$

where $B^{\text{AMGS}}$ is a base signifying the asymptotic convergence rate and $C^{\text{AMGS}}$ is a proportionality constant, given by

$$B^{\text{AMGS}} \stackrel{\text{def}}{=} \left(1 + \sqrt{\frac{\mu}{2(L_f - \mu)}}\right)^2 = \left(1 + \sqrt{\frac{q}{2(1-q)}}\right)^2,$$
$$C^{\text{AMGS}} \stackrel{\text{def}}{=} \frac{1}{\left(\sqrt{L_f - \mu} + \sqrt{\frac{\mu}{2}}\right)^2} = \frac{1}{L_f\left(\sqrt{1-q} + \sqrt{\frac{q}{2}}\right)^2}.$$

For FGM, with the same parameter choice as in the non-strongly convex case, we have

$$A_k^{\text{FGM}} = A_i^{\text{FGM}} \geq C^{\text{FGM}}\left(B^{\text{FGM}}\right)^i,$$

where

$$B^{\text{FGM}} \stackrel{\text{def}}{=} \left(\frac{1}{1-\sqrt{q}}\right)^2, \qquad C^{\text{FGM}} \stackrel{\text{def}}{=} \frac{1}{2L_f}.$$

The relation

$$\sqrt{\frac{B^{\text{AMGS}}}{B^{\text{FGM}}}} = \frac{1 + \sqrt{\frac{q}{2(1-q)}}}{\frac{1}{1-\sqrt{q}}} = \frac{1 - q + \sqrt{\frac{q(1-q)}{2}}}{1 + \sqrt{q}} < \frac{1 + \sqrt{\frac{q}{2}}}{1 + \sqrt{q}} < 1$$

shows that FGM is asymptotically more efficient, by a considerable margin, than AMGS. The reason for this is that lower bound at iteration $k$, $w_{k+1}(\boldsymbol{x})$, employed by FGM is an accurate approximation of the objective function around a recently used point. For non-strongly convex objectives, AMGS has a lower bound that is tangent at the current iterate, which is slightly more accurate to that of FGM but requires 2 WTU to update, giving it no computational advantage. The AMGS bound is very poor in the strongly convex case because the lower bound square term $\frac{\mu}{2}\|\boldsymbol{x} - \boldsymbol{x_0}\|_2^2$ is centered around $\boldsymbol{x}_0$ and does not take into account recent information.

**Appendix C. Proof of Proposition 6.**

With judicious use of parameters $r_u$ and $r_d$, the average WTU cost of an iteration of ACGM can be adjusted to equal that of FISTA (also evidenced in subsection 4.1). Consequently, it is adequate to compare the convergence guarantees of the two algorithms when indexed in iterations.

Combining (86) and (82), we obtain

$$A_{k+1}^{\text{ACGM}} = \left(\sqrt{\frac{1}{4L_{k+1}^{\text{ACGM}}}} + \sqrt{\frac{1}{4L_{k+1}^{\text{ACGM}}} + A_k^{\text{ACGM}}}\right)^2.$$

Replacing (86) in ACGM with

$$t_{k+1} = \frac{1 + \sqrt{1 + 4t_k^2}}{2}, \quad k \geq 0, \tag{93}$$

results in an algorithm that produces identical iterates to FISTA. The convergence analysis of ACGM yields the corresponding accumulated weight for FISTA

$$A_{k+1}^{\text{FISTA}} = \left(\sqrt{\frac{1}{4L_{k+1}^{\text{FISTA}}}} + \sqrt{\frac{1}{4L_{k+1}^{\text{FISTA}}} + \frac{L_k^{\text{FISTA}}}{L_{k+1}^{\text{FISTA}}} A_k^{\text{FISTA}}}\right)^2.$$

Both methods start with the same state $A_0^{\text{ACGM}} = A_0^{\text{FISTA}} = 0$. The line-search procedure of ACGM is guaranteed to produce Lipschitz constant estimates no greater than those of FISTA for the same local curvature, i.e., $L_k^{\text{ACGM}} \leq L_k^{\text{FISTA}}$, $k \geq 0$. FISTA, by design, can only accommodate a Lipschitz constant estimate increase, namely $L_k^{\text{FISTA}} \leq L_{k+1}^{\text{FISTA}}$, $k \geq 0$. Thus, regardless of the fluctuation in the local curvature of $f$, we have

$$A_k^{\text{ACGM}} \geq A_k^{\text{FISTA}}, \quad k \geq 0.$$